\documentclass[11pt]{amsart}
\usepackage{amssymb,verbatim}
\usepackage{amsthm}
\usepackage[mathscr]{eucal}
\usepackage[OT1]{fontenc}
\usepackage{url}

\usepackage{pstricks,pst-node,pst-text,pst-3d}
\usepackage{gastex}
\usepackage{epsf}

\newtheorem{thm}{Theorem}[section]
\newtheorem{lem}[thm]{Lemma}
\newtheorem{cor}[thm]{Corollary}
\newtheorem{pro}[thm]{Proposition}

\theoremstyle{remark}
\newtheorem{eg}[thm]{Example}
\newtheorem{notn}[thm]{Notation}

\newenvironment{proof*}{\trivlist\item[\hskip\labelsep\emph{Proof.}]}{\endtrivlist}

\DeclareSymbolFont{rsfscript}{OMS}{rsfs}{m}{n}
\DeclareSymbolFontAlphabet{\mathrsfs}{rsfscript}

\newcommand{\Cg}{\mathsf{C}}

\newcommand{\Gg}{\mathsf{G}}
\newcommand{\Hg}{\mathsf{H}}
\newcommand{\Kg}{\mathsf{K}}

\newcommand{\Pg}{\mathsf{P}}

\newcommand{\Rg}{\mathsf{R}}
\newcommand{\Sg}{\mathsf{S}}

\newcommand{\FAAR}{\mathrm{FAAR}}
\newcommand{\SAAR}{\mathrm{SAAR}}
\newcommand{\fapprox}{\stackrel\FAAR\approx}
\newcommand{\sapprox}{\stackrel\SAAR\approx}

\newcommand{\Lgr}{\mathbin{\mathrsfs{L}}}
\newcommand{\Rgr}{\mathbin{\mathrsfs{R}}}
\newcommand{\Hgr}{\mathbin{\mathrsfs{H}}}

\newcommand{\A}{\mathrm{A}}

\newcommand{\Bigand}{\operatornamewithlimits{\hbox{\LARGE$\&$}}}
\newcommand{\bigand}{\operatornamewithlimits{\hbox{\large$\&$}}}

\newcommand{\Ad}{\operatorname{A}}
\newcommand{\Syn}{\operatorname{Syn}}

\newcommand{\up}[1]{\textup{#1}}

\begin{document}

\title[The algebra of adjacency patterns]{The algebra of adjacency patterns:\\
Rees matrix semigroups with reversion}

\author{Marcel Jackson}
\address{La Trobe University, Victoria 3086,  Australia}
\email{M.G.Jackson@latrobe.edu.au}
\author{Mikhail Volkov}
\address{Ural State University, Ekaterinburg 620083, Russia}
\email{Mikhail.Volkov@usu.ru}


\thanks{The first author was supported by ARC Discovery Project Grant DP0342459.
The second author was supported by the Program 2.1.1/3537 of the
Russian Education Agency}

\begin{abstract}
   We establish a surprisingly close relationship between universal Horn
   classes of directed graphs and varieties generated by so-called adjacency
   semigroups which are Rees matrix semigroups over the trivial group with
   the unary operation of reversion. In particular, the lattice of subvarieties
   of the variety generated by adjacency semigroups that are regular unary
   semigroups is essentially the same as the lattice of universal Horn classes
   of reflexive directed graphs. A number of examples follow, including a limit
   variety of regular unary semigroups and finite unary semigroups with
   \textsf{NP}-hard variety membership problems.
\end{abstract}

\maketitle

\section*{Introduction and overview}

The aim of this paper is to establish and to explore a new link between
graph theory and algebra. Since graphs form a universal language of discrete
mathematics, the idea to relate graphs and algebras appears to be natural,
and several useful links of this kind can be found in the literature.
We mean, for instance, the graph algebras of McNulty and Shallon \cite{mcnsha},
the closely related flat graph algebras \cite{will96}, and ``almost trivial'' algebras
investigated in \cite{jezand others,jezmck} amongst other places. While each
of the approaches just mentioned has proved to be useful and has yielded
interesting applications, none of them seem to share two important features
of the present contribution. The two features can be called naturalness and
surjectivity.

Speaking about naturalness, we want to stress that the algebraic objects (adjacency
semigroups) that we use here to interpret graphs have not been invented for this
specific purpose. Indeed, adjacency semigroups belong to a well established
class of unary semigroups\footnote{Here and below the somewhat oxymoronic term
``unary semigroup'' abbreviates the precise but longer expression ``semigroup
endowed with an extra unary operation''.} that have been considered by many
authors. We shall demonstrate how graph theory both sheds a new light on
some previously known algebraic results and provides their extensions and
generalizations. By surjectivity we mean that, on the level of appropriate
classes of graphs and unary semigroups, the interpretation map introduced
in this paper becomes ``nearly'' onto; moreover, the map induces a lattice
isomorphism between the lattices of such classes provided one excludes just
one element on the semigroup side. This implies that our approach allows one
to interpret both graphs within unary semigroups and unary semigroups within
graphs.

The paper is structured as follows. In Section~\ref{sec:graphs} we recall some
notions related to graphs and their classes and present a few results and examples
from graph theory that are used in the sequel. Section~\ref{sec:main} contains our
construction and the formulations of our main results: Theorems~\ref{thm:fund}
and~\ref{thm:fundref}. These theorems are proved in Sections~3 and~4 respectively
while Section~5 collects some of their applications.

We assume the reader's acquaintance with basic concepts of universal algebra
and first-order logics such as ultraproducts or the HSP-theorem, see, e.g.,
\cite{BuSa81}. As far as graphs and semigroups are concerned, we have tried
to keep the presentation to a reasonable extent self-contained. We do occasionally
mention some non-trivial facts of semigroup theory but only in order to place
our considerations in a proper perspective. Thus, most of the material should
be accessible to readers with very basic semigroup-theoretic background
(such as some knowledge of Green's relations and of the Rees matrix
construction over the trivial group, cf.~\cite{how}).

\section{Graphs and their classes}
\label{sec:graphs}

In this paper, \emph{graph} is a structure $\Gg:=\langle V;\sim \rangle$,
where $V$ is a set and $\sim\ \subseteq V\times V$ is a binary relation.
In other words, \textbf{we consider all graphs to be directed}, and do not allow
multiple edges (but do allow loops). Of course, $V$ is often referred to
as the set of \emph{vertices} of the graph and $\sim$ as the set of \emph{edges}.
As is usual, we write $a\sim b$ in place of $(a,b)\in {\sim}$.  Conventional
undirected graphs are essentially the same as graphs whose edge relation
is symmetric (satisfying $x\sim y\rightarrow y\sim x$), while a \emph{simple
graph} is a symmetric graph without loops. It is convenient for us to allow
the empty graph $\underline{\bf 0}:=\langle\varnothing;\varnothing\rangle$.

All classes of graphs that come to consideration in this paper are \emph{universal
Horn classes}. We recall their definition and some basic properties. Of course,
the majority of the statements below are true for arbitrary structures, but our
interest is only in the graph case.  See Gorbunov~\cite{gor} for more details.

Universal Horn classes can be defined both syntactically (via specifying an
appropriate sort of first order formulas) and semantically (via certain class
operators). We first introduce the operator definition for which we recall
notation for a few standard class operators. The operator for taking isomorphic
copies is $\mathbb{I}$. We use $\mathbb{S}$ to denote the operator taking
a class $K$ to the class of all substructures of structures in $K$; in the case
when $K$ is a class of graphs, substructures are just induced subgraphs of graphs
in $K$. Observe that the empty graph $\underline{\bf 0}$ is an induced subgraph
of any graph and thus belongs to any $\mathbb{S}$-closed class of graphs. We denote
by $\mathbb{P}$ the operator of taking direct products. For graphs, we allow the
notion of an empty direct product, which we identify (as is the standard convention)
with the 1-vertex looped graph $\underline{\bf 1}:=\langle \{0\};\{(0,0)\}\rangle$.
If we exclude the empty product, we obtain the operator $\mathbb{P}^{+}$ of taking
nonempty direct products. By $\mathbb{P}_{\mathrm{u}}$ we denote the operator
of taking ultraproducts. Note that ultraproducts---unlike direct products---are
automatically nonempty.

A class $K$ of graphs is an \emph{universal Horn class} if $K$ is closed
under each of the operators $\mathbb{I}$, $\mathbb{S}$, $\mathbb{P}^{+}$,
and $\mathbb{P}_{\mathrm{u}}$. In the sequel, we write ``uH class'' in
place of ``universal Horn class''. It is well known that the least uH class
containing a class $L$ of graphs is the class $\mathbb{ISP^{+}P}_{\mathrm{u}}(L)$
of all isomorphic copies of induced subgraphs of nonempty direct products of
ultraproducts of $L$; this uH class is referred to as the uH class
\emph{generated} by $L$.

If the operator $\mathbb{P}^{+}$ in the above definition is extended to $\mathbb{P}$,
then one obtains the definition of a \emph{quasivariety} of graphs. The quasivariety
\emph{generated} by a given class $L$ is known to be equal to $\mathbb{ISPP}_{\mathrm{u}}(L)$.
It is not hard to see that $\mathbb{ISPP}_{\mathrm{u}}(L)=
\mathbb{I}(\mathbb{ISP^{+}P}_{\mathrm{u}}(L)\cup\{\underline{\bf 1}\})$, showing that
there is little or no difference between the uH class and the quasivariety generated
by $L$. However, as examples described later demonstrate, there are many well studied
classes of graphs that are uH classes but not quasivarieties.

As mentioned, uH classes also admit a well known syntactic characterization.
An \emph{atomic formula} in the language of graphs is an expression of the form
$x\sim y$ or $x\approx y$ (where $x$ and $y$ are possibly identical variables).
A \emph{universal Horn sentence} (abbreviated to ``uH sentence'') in the language
of graphs is a sentence of one of the following two forms (for some
$n\in\omega:= \{0, 1, 2, \dots\}$):
\[
(\forall x_{1}\forall x_{2}\ldots)\left(\left(\Bigand_{1\leq i\leq
n}\Phi_{i}\right)\rightarrow \Phi_{0}\right)\quad \mbox{ or }\quad
(\forall x_{1}\forall x_{2}\ldots)\left(\bigvee_{0\leq i\leq
n}\neg\Phi_{i}\right)
\]
where the $\Phi_{i}$ are atomic, and $x_{1},x_{2},\ldots$ is a list
of all variables appearing.  In the case when $n=0$, a uH sentence
of the first kind is simply the universally quantified atomic expression
$\Phi_{0}$.  Sentences of the first kind are usually called \emph{quasi-identities}.
As is standard, we omit the universal quantifiers when describing uH sentences;
also the expressions $x\not\approx y$ and $x\nsim y$ abbreviate $\neg x\approx y$
and $\neg x\sim y$ respectively. Satisfaction of uH sentences by graphs is defined
in the obvious way. We write $\Gg\models\Phi$ ($K\models\Phi$) to denote
that the graph $\Gg$ (respectively, each graph in the class $K$) satisfies
the uH sentence~$\Phi$.

The Birkhoff theorem identifying varieties of algebras with equationally
defined classes has a natural analogue for uH classes, which is usually
attributed to Mal'cev.  Here we state it in the graph setting.
\begin{lem}\label{lem:malcev}
    A class $K$ of graphs is a uH class if and only if it is the class
    of all models of some set of uH sentences.
\end{lem}
In particular, the uH class $\mathbb{ISP^{+}P}_\mathrm{u}(L)$ generated by a class
$L$ is equal to the class of models of the uH sentences holding in $L$.

Recall that we allow the empty graph $\underline{\bf 0}:=\langle \varnothing;\varnothing\rangle$.
Because there are no possible variable assignments into the empty set, $\underline{\bf 0}$
can fail no uH sentence and hence lies in every uH class. Thus, allowing $\underline{\bf 0}$
brings the advantage that the collection of all uH classes forms a lattice whose meet is intersection:
$A\wedge B:=A\cap B$ and whose join is given by $A\vee B:=\mathbb{ISP^{+}P}_{\mathrm{u}}(A\cup B)$.
Furthermore, the inclusion of $\underline{\bf 0}$ allows every set of uH sentences to have a model
(for example, the contradiction $x\not\approx x$ axiomatizes the class $\{\underline{\bf 0}\}$).
In the world of varieties of algebras, it is the one element algebra that plays these roles.

When $\mathbb{IP}_{\mathrm{u}}(L)=\mathbb{I}(L)$ (such as when $L$ consists of finitely many
finite graphs), we have $\mathbb{ISP^{+}P}_{\mathrm{u}}(L)=\mathbb{ISP^{+}}(L)$, and there is
a handy structural characterization of the uH class generated by $L$.

\begin{lem}\label{lem:sep}
    Let $L$ be an ultraproduct closed class of graphs and let
    $\Gg$ be a graph.  We have $\Gg\in\mathbb{ISP^{+}P}_{\mathrm{u}}(L)$ if and only
    if there is at least one homomorphism from $\Gg$ into a member of $L$ and the following
    two     \emph{separation conditions} hold\up:
    \begin{enumerate}
    \item for each pair of distinct vertices $a,b$ of $G$, there
    is $\Hg\in L$ and a homomorphism  $\phi:\Gg\to \Hg$ with $\phi(a)\neq
    \phi(b)$\up;
    \item for each pair of vertices $a,b$ of $G$ with $a\nsim b$ in $\Gg$,
    there is  $\Hg\in L$ and a homomorphism $\phi:\Gg\to \Hg$ with $\phi(a)\nsim\phi(b)$ in $\Hg$.
    \end{enumerate}
\end{lem}

The 1-vertex looped graph $\underline{\bf 1}$ always satisfies the two separation conditions,
yet it fails every uH sentence of the second kind; this is why the lemma asks additionally that
there be at least one homomorphism from $\Gg$ into some member of $L$.  If $\Gg=\underline{\bf 1}$
and no such homomorphism exists, then evidently, each member of $L$ has nonlooped vertices, and so
$L\models x\nsim x$, a law failing on $\underline{\bf 1}$.  Hence
$\underline{\bf1}\notin\mathbb{ISP^{+}P}_{\mathrm{u}}(L)$ by Lemma~\ref{lem:malcev}.
Conversely, if there is such a homomorphism, then $\underline{\bf 1}$ is isomorphic to an
induced subgraph of some member of $L$ and hence $\underline{\bf1}\in\mathbb{ISP^{+}P}_{\mathrm{u}}(L)$.
If the condition that there is at least one homomorphism from $\Gg$ into some member of $L$ is
dropped, then Lemma \ref{lem:sep} instead characterizes membership in the quasivariety generated by $L$.

We now list some familiar uH sentences.
\begin{itemize}
    \item reflexivity: $x\sim x$,
    \item anti-reflexivity: $x\not\sim x$,
    \item symmetry: $x\sim y\rightarrow y\sim x$,
    \item anti-symmetry: $x\sim y\And y\sim x\rightarrow x\approx y$,
    \item transitivity: $x\sim y\And y\sim z\rightarrow x\sim y$.
\end{itemize}
All except anti-reflexivity are quasi-identities.

These laws appear in many commonly investigated classes of graphs.
We list a number of examples that are of interest later in the
paper (mainly in its application part, see Section~\ref{sec:applications}
below).

 \begin{eg} \label{eg:preorders}
     Preorders.
    \end{eg}

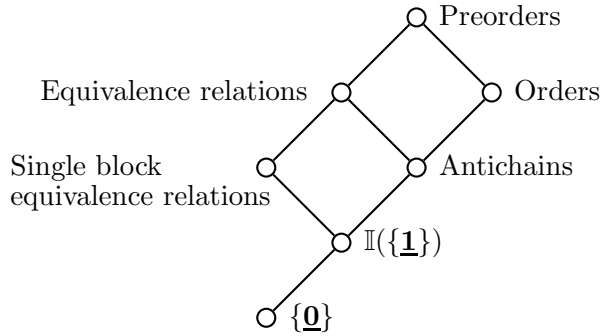
\begin{figure}[ht]
\begin{center}
\begin{pspicture}(3,4)(0,0)
      \cnodeput(0,0){a}{}
      \cnodeput(1,1){b}{}
      \cnodeput(0,2){c}{}
      \cnodeput(2,2){d}{}
      \cnodeput(1,3){e}{}
      \cnodeput(3,3){f}{}
      \cnodeput(2,4){g}{}
      \psset{arrows=-,nodesep=0}
      \ncline{a}{b}
      \ncline{b}{c}
      \ncline{b}{d}
      \ncline{c}{e}
      \ncline{e}{g}
      \ncline{d}{f}
      \ncline{d}{e}
      \ncline{f}{g}
      \put(0.3,-0.1){$\{\underline{\bf 0}\}$}
      \put(1.3,0.9){$\mathbb{I}(\{\underline{\bf 1}\})$}
      \put(-3.4,1.9){Single block}
      \put(-3.4,1.5){equivalence relations}
      \put(2.3,1.9){Antichains}
      \put(-3,2.9){Equivalence relations}
      \put(3.3,2.9){Orders}
      \put(2.3,3.9){Preorders}
\end{pspicture}
\end{center}
\caption{The lattice of uH classes of preorders}\label{pic:preorders}
\end{figure}

    This class is defined by reflexivity and
    transitivity and is a quasivariety.  Some
    well known subclasses are:
    \begin{itemize}
        \item equivalence relations (obtained by adjoining the
        symmetry law);
        \item partial orders (obtained by adjoining the
        anti-symmetry law);
        \item anti-chains (the intersection of partial orders and
        equivalence relations);
        \item complete looped graphs, or equivalently, single
        block equivalence relations (axiomatized by $x\sim y$).
    \end{itemize}
    In fact it is easy to see that, along with the 1-vertex partial
    orders and the trivial class $\{\underline{\bf 0}\}$, this exhausts
    the list of all uH classes of preorders, see Fig.\,\ref{pic:preorders}
    (the easy proof is sketched before Corollary 6.4 of \cite{CDJP}, for
    example).

      \begin{eg}\label{eg:simple}
      Simple (that is, anti-reflexive and symmetric) graphs.
      \end{eg}

    Sub-uH~classes of simple graphs have been heavily investigated, and include some
    very interesting families. In order to describe some of these families, we need
    a sequence of graphs introduced by Ne\v{s}et\v{r}il and Pultr~\cite{nespul}. For
    each integer $k\ge 2$, let $\Cg_{k}$ denote the graph on the vertices
    $0,\ldots,k+1$ obtained from the complete loopless graph on these vertices
    by deleting the edges (in both directions) connecting $0$ and $k+1$, $0$ and $k$,
    and $1$ and $k+1$.
\begin{figure}[hb]
\begin{center}
\begin{pspicture}(8.5,2.3)(-6,0)
      \cnodeput(-4,2){b1}{$0$}
      \cnodeput(-2,2){c1}{$2$}
      \cnodeput(-3,0){d1}{$1$}
      \cnodeput(-1,0){a1}{$3$}
      \cnodeput(1,2){b}{$0$}
      \cnodeput(3,2){c}{$2$}
      \cnodeput(5,2){d}{$4$}
      \cnodeput(2,0){e}{$1$}
      \cnodeput(4,0){a}{$3$}
      \psset{arrows=-,nodesep=0}
      \ncline{b1}{d1}
      \ncline{c1}{d1}
      \ncline{c1}{a1}
      \ncline{a}{d}
      \ncline{b}{c}
      \ncline{c}{d}
      \ncline{b}{e}
      \ncline{c}{a}
      \ncline{e}{a}
      \ncline{c}{e}
\end{pspicture}
\end{center}
\caption{Graphs $\Cg_{2}$ and $\Cg_{3}$}\label{pic:c2 and c3}
\end{figure}
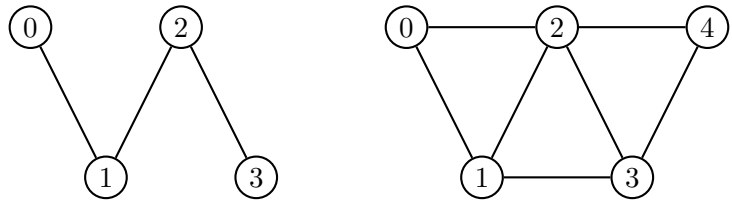
Fig.\,\ref{pic:c2 and c3} shows the graphs $\Cg_{2}$ and $\Cg_{3}$; here and
below we adopt the convention that an undirected edge between two vertices, say $a$
and $b$, represents two directed edges $a\sim b$ and $b\sim a$.

Recall that a simple graph $\Gg$ is said to be $n$-\emph{colorable} is there exists
a homomorphism from $\Gg$ into the complete loopless graph on $n$ vertices.

    \begin{eg}\label{eg:2-color}
    The $2$-colorable graphs (equivalently, bipartite graphs).
    \end{eg}

    This is the uH class $\mathbb{ISP}^{+}(\Cg_{2})$ generated by the graph $\Cg_{2}$
    (Ne\v{s}et\v{r}il and Pultr \cite{nespul}) and has no finite axiomatization. Caicedo
    \cite{cai} showed that the lattice of sub-uH classes of $\mathbb{ISP}^{+}(\Cg_{2})$
    is a 6-element chain:  besides $\mathbb{ISP}^{+}(\Cg_{2})$, it contains the class
    of disjoint unions of complete bipartite graphs, which is axiomatized within simple
    graphs by the law
    $$x_{0}\sim x_{1}\And x_{1}\sim x_{2}\And x_{2}\sim
    x_{3}\rightarrow x_{0}\sim x_{3};$$ the class of disjoint
    unions of paths of length at most $1$ (axiomatized within simple graphs by
    $x\nsim y\vee y\nsim z$); the edgeless graphs (axiomatized by $x\nsim y$),
    the 1-vertex edgeless graphs ($x\approx y$); and the trivial
    class $\{\underline{\bf 0}\}$.

    Every finite simple graph either lies in a sub-uH class of $\mathbb{ISP}^{+}(\Cg_{2})$
    or generates a uH class that: 1) is not finitely axiomatizable, 2) contains $\mathbb{ISP}^{+}(\Cg_{2})$,
    and 3) has uncountably many sub-uH classes \cite[Theorem 4.7]{gorkra}, see also~\cite{kra}.

    \begin{eg}\label{eg:k-color}
    The $k$-colorable graphs.
    \end{eg}

    More generally, Nesetril and Pultr \cite{nespul} showed that for any
    $k\geq 2$, the class of all $k$-colorable graphs is the uH class generated
    by $\Cg_{k}$. These classes have no finite basis for their uH
    sentences and for $k>2$ have \textsf{NP}-complete finite membership
    problem, see \cite{garey}.

    \begin{eg}\label{eg:graphgen}
    A generator for the class $\underline{\mathsf{G}}$ of all graphs.
    \end{eg}

    The class of all graphs is generated as a uH class by a single finite graph. Indeed,
    it is trivial to see that for any graph $\Gg$, there is a family of 3-vertex graphs
    such that the separation conditions of Lemma~\ref{lem:sep} hold. Since there are only
    finitely many non-isomorphic 3-vertex graphs, any graph containing these as induced
    subgraphs generates the uH class of all graphs. Alternatively, the reader can easily
    verify using Lemma~\ref{lem:sep}  that the following graph $\Gg_{1}$ generates
    the uH class of all graphs:

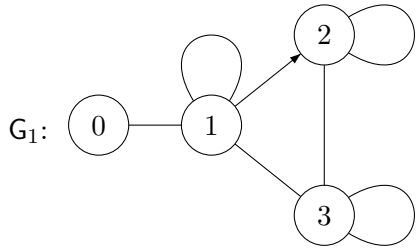
\begin{figure}[hb]
\unitlength=.6mm
\begin{center}
\begin{picture}(60,45)
\node(B)(0,20){0}
\node(C)(25,20){1}
\node(D)(50,40){2}
\node(E)(50,0){3}
\drawedge(C,D){}
\gasset{AHnb=0}
\drawloop(C){}
\drawloop[loopangle=0](D){}
\drawloop[loopangle=0](E){}
\drawedge(B,C){}
\drawedge(D,E){}
\drawedge(C,E){}
\put(-20,17){$\Gg_{1}$:}
\end{picture}
\end{center}
\caption{Generator for the uH class of all graphs}\label{pic:gen_for_all}
\end{figure}

    \begin{eg}\label{eg:symmetric}
    A generator for the class $\underline{\mathsf{G}}_{\mathrm{symm}}$ of all symmetric graphs.
    \end{eg}

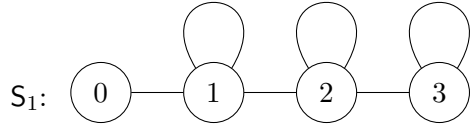
\begin{figure}[ht]
\unitlength=.6mm
\begin{center}
\begin{picture}(75,15)
\node(B)(0,0){0}
\node(C)(25,0){1}
\node(D)(50,0){2}
\node(E)(75,0){3}
\gasset{AHnb=0}
\drawloop(C){}
\drawloop(D){}
\drawloop(E){}
\drawedge(B,C){}
\drawedge(C,D){}
\drawedge(D,E){}
\put(-20,-3){$\Sg_{1}$:}
\end{picture}
\end{center}
\caption{Generator for the uH class of all symmetric graphs}\label{pic:gen_for_sym}
\end{figure}

    Using Lemma \ref{lem:sep}, it is easy to prove that the class of symmetric graphs
    is generated as a uH class by the graph $\Sg_{1}$ shown in Fig.\,\ref{pic:gen_for_sym}.

    \begin{eg}\label{eg:simplegraphgen}
    The class of simple graphs has no finite generator.
    \end{eg}

    The class of all simple graphs is not generated by any
    finite graph, since a finite graph on $n$ vertices is
    $n$-colorable, while for every positive integer $n$ there is a simple graph
    that is not $n$-colorable (the complete simple graph on
    $n+1$ vertices, for example). However the uH class generated
    by the following 2-vertex graph $\Sg_{2}$ contains all simple graphs
    (this is well known and follows easily using Lemma \ref{lem:sep}).

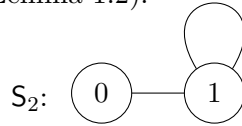
\begin{figure}[ht]
\unitlength=.6mm
\begin{center}
\begin{picture}(25,12)(0,-2)
\node(B)(0,0){0}
\node(C)(25,0){1}
\gasset{AHnb=0}
\drawloop(C){}
\drawedge(B,C){}
\put(-20,-3){$\Sg_{2}$:}
\end{picture}
\end{center}
\caption{2-vertex graph whose uH class contains all simple graphs}\label{pic:gen_for_sim}
\end{figure}

    \begin{eg}\label{eg:refgen}
    A generator for the class $\underline{\mathsf{G}}_{\mathrm{ref}}$ of all reflexive graphs.
    \end{eg}

    The class of reflexive graphs is generated by the following graph $\Rg_{1}$,
    while the class of reflexive and symmetric graphs is generated by the graph
    $\mathsf{RS}_{1}$.

\begin{figure}[hb]
\unitlength=.6mm
\begin{center}
\begin{picture}(175,40)(-2,0)
\node(C)(25,20){0}
\node(D)(50,38){1}
\node(E)(50,2){2}
\drawedge(C,D){}
\gasset{AHnb=0}
\drawloop[loopangle=180](C){}
\drawloop[loopangle=0](D){}
\drawloop[loopangle=0](E){}
\drawedge(D,E){}
\drawedge(C,E){}
\put(-8,17){$\Rg_{1}$:}
\node(C1)(125,20){0}
\node(D1)(150,20){1}
\node(E1)(175,20){2}
\drawloop(C1){}
\drawloop(D1){}
\drawloop(E1){}
\drawedge(C1,D1){}
\drawedge(D1,E1){}
\put(100,17){$\mathsf{RS}_{1}$:}
\end{picture}
\end{center}
\caption{Generators for reflexive and reflexive-symmetric graphs}\label{pic:gen_for_ref}
\end{figure}
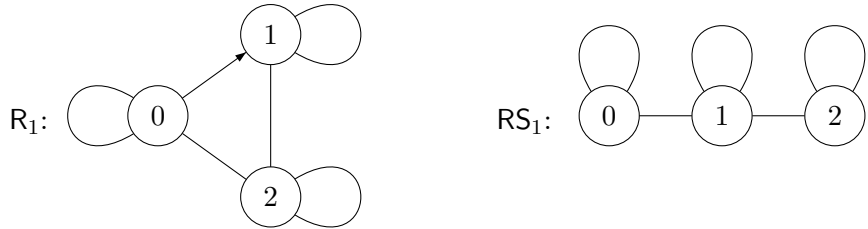

\section{The adjacency semigroup of a graph}
\label{sec:main}

Given a graph $\Gg=\langle V;\sim\rangle$, its \emph{adjacency semigroup}
$\A(\Gg)$ is defined on the set $(V\times V)\cup\{0\}$ and
the multiplication rule is
\begin{gather*}
(x,y)(z,t)=\begin{cases}
(x,t)&\ \text{if}\ y\sim z,\\
0    &\ \text{if}\ y\nsim z;
\end{cases}\\
a0=0a=0\ \text{ for all $a\in\A(\Gg)$}.
\end{gather*}
In terms of semigroup theory, $\A(\Gg)$ is the \emph{Rees matrix semigroup
over the trivial group} using the adjacency matrix of the graph $\Gg$
as a sandwich matrix. We describe here the Rees matrix construction
in a specific form that is used in the present paper.

Let $I,J$ be nonempty sets and $0\notin I\cup J$. Let $P=(P_{i,j})$ be
a $J\times I$ matrix (the \emph{sandwich matrix}) over the set $\{0,1\}$.
The \emph{Rees matrix semigroup over the trivial group} $M^0[P]$
is the semigroup on the set $(I\times J)\cup \{0\}$ with multiplication
\begin{gather*}
a\cdot 0=0\cdot a=0\ \text{ for all $a\in (I\times J)\cup \{0\}$}, \text{ and}\\
(i_1,j_1)\cdot(i_2,j_2)=\begin{cases}0&\mbox{ if }P_{j_1,i_2}=0,\\
(i_1,j_2)&\mbox{ if }P_{j_1,i_2}=1.
\end{cases}
\end{gather*}
The Rees-Sushkevich Theorem (see \cite[Theorem 3.3.1]{how})
states that, up to isomorphism, the completely 0-simple semigroups
with trivial subgroups are precisely the Rees matrix semigroups
over the trivial group and for which each row and each column of the
sandwich matrix contains a nonzero element. If the matrix $P$ has no
0 entries, then the set $M[P]=M^0[P]\setminus\{0\}$ is a subsemigroup.
Semigroups of the form $M[P]$ are called \emph{rectangular bands},
and they are precisely the completely simple semigroups with trivial
subgroups.

Back to adjacency semigroups, we always think of $\A(\Gg)$ as endowed
with an additional unary operation $a\mapsto a'$ which we call \emph{reversion}
and define as follows:
$$(x,y)'=(y,x), \quad 0'=0.$$
Notice that by this definition $(a')'=a$ for all $a\in\A(\Gg)$.

The main contribution in this paper is the fact that uH classes of graphs
are in extremely close correspondence with unary semigroup varieties generated
by adjacency semigroups, and our proof of this will involve a translation of
uH~sentences of graphs into unary semigroup identities. However, before we
proceed with precise formulations and proofs of general results, the reader
may find it useful to check that several of the basic uH sentences used
in Section~\ref{sec:graphs} correspond via the adjacency semigroup construction
to rather natural semigroup-theoretic properties. Indeed, all the following are
quite easy to verify:
\begin{itemize}
    \item reflexivity of $\Gg$ is equivalent to $\A(\Gg)\models xx'x\approx x$;
    \item anti-reflexivity of $\Gg$ is equivalent to $\A(\Gg)\models
    xx'z\approx zxx'\approx xx'$ (these laws can be abbreviated to
    $xx'\approx 0$);
    \item symmetry of $\Gg$ is equivalent to $\A(\Gg)\models (xy)'\approx y'x'$;
    \item $\Gg$ is empty (satisfies $x\not\approx x$) if and only if $\A(\Gg)\models
    x\approx y$;
    \item $\Gg$ has one vertex (satisfies $x\approx y$) if and only if $\A(\Gg)\models
    x\approx x'$;  also, $\Gg$ is the one vertex looped graph (satisfies $x\sim y$) if and
    only if $\A(\Gg)$ additionally satisfies $xx\approx x$.
\end{itemize}

Observe that the unary semigroup identities that appear in the above
examples are in fact used to define the most widely studied types of
semigroups endowed with an extra unary operation modelling various
notions of an inverse in groups. For instance, a semigroup satisfying
the identities
\begin{equation}
\label{eq:repeat}
x''\approx x
\end{equation}
(which always holds true in adjacency semigroups) and
\begin{equation}
\label{eq:involution}
(xy)'\approx y'x'
\end{equation}
(which is a semigroup counterpart of symmetry) is called \emph{involution
semigroup} or \emph{$*$-semigroup}. If such a semigroup satisfies also
\begin{equation}
\label{eq:regular}
xx'x\approx x
\end{equation}
(which corresponds to reflexivity), it is called a \emph{regular $*$-semigroup}.
Semigroups satisfying \eqref{eq:repeat} and \eqref{eq:regular} are called
\emph{I-semigroups} in Howie \cite{how}; note that an I-semigroup satisfies
$x'xx'\approx x'x''x'\approx x'$, so that $x'$ is an inverse of $x$.
Semigroups satisfying \eqref{eq:regular} are often called \emph{regular
unary semigroups}. There exists vast literature on all these types of
unary semigroups; clearly, the present paper is not a proper place to
survey this literature but we just want to stress once more that the
range of the adjacency semigroup construction is no less natural than
its domain.

When $K$ is a class of graphs, we use the notation $\A(K)$ to denote
the class of all adjacency semigroups of members of $K$. As usual,
the operator of taking homomorphic images is denoted by $\mathbb{H}$.
We let $\mathscr{A}$ denote the variety $\mathbb{HSP}(\A(\underline{\mathsf{G}}))$
generated by all adjacency semigroups of graphs, and let
$\mathscr{A}_\mathrm{ref}$ and $\mathscr{A}_\mathrm{symm}$
denote the varieties $\mathbb{HSP}(\A(\underline{\mathsf{G}}_\mathrm{ref}))$ and
$\mathbb{HSP}(\A(\underline{\mathsf{G}}_\mathrm{symm}))$ generated by all adjacency semigroups
of reflexive graphs and of symmetric graphs respectively.

Our first main result is:
\begin{thm}\label{thm:fund}
Let $K$ be any nonempty class of graphs and let $\Gg$ be a graph.
The graph $\Gg$ belongs to the uH-class generated by $K$ if and only
if the adjacency semigroup $\A(\Gg)$ belongs to the variety generated by
the adjacency  semigroups $\A(\Hg)$ with $\Hg\in K$.
\end{thm}

This immediately implies that the assignment $\Gg\mapsto \A(\Gg)$ induces
an injective join-preserving map from the lattice of all uH-classes of graphs
to the subvariety lattice of the variety $\mathscr{A}$. The latter fact can
be essentially refined for the case of reflexive graphs. In order to describe
this refinement, we need an extra definition.

Let $I$ be a nonempty set. We endow the set $B=I\times I$ with a unary semigroup
structure whose multiplication is defined by
$$(i,j)(k,\ell)=(i,\ell)$$
and whose unary operation is defined by
$$(i,j)'=(j,i).$$
It is easy to check that $B$ becomes a regular $*$-semigroup. We call regular
$*$-semigroups constructed this way \emph{square bands}. Clearly, square bands
satisfy
\begin{equation}
\label{eq:square bands}
x^2\approx x\ \text{ and }\ xyz\approx xz,
\end{equation}
and in fact it can be shown that the class $\mathscr{SB}$ of all square bands
constitutes a variety of unary semigroups defined within the variety of all
regular $*$-semigroups by the identities \eqref{eq:square bands}.

Let $L(\underline{\mathsf{G}}_\mathrm{ref})$ denote the lattice of
sub-uH classes of $\underline{\mathsf{G}}_\mathrm{ref}$ and let
$L(\mathscr{A}_\mathrm{ref})$ denote the lattice of subvarieties
of $\mathscr{A}_\mathrm{ref}$. Let $L^{+}$ denote the result of
adjoining a new element $\mathsf{\underline{S}}$ to
$L(\underline{\mathsf{G}}_\mathrm{ref})$ between the class of
single block equivalence relations and the class containing
the empty graph. (The reader may wish to look at Fig.~\ref{pic:preorders}
to see the relative location of these two uH classes.)
Meets and joins are extended to $L^{+}$ in the weakest way.
So $L^{+}$ is a lattice in which $L(\underline{\mathsf{G}}_\mathrm{ref})$
is a sublattice containing all but one element.

We are now in a position to formulate our second main result.
\begin{thm}\label{thm:fundref}
    Let $\iota$ be the map from $L^{+}$ to $L(\mathscr{A}_\mathrm{ref})$
    defined by $\mathsf{\underline{S}}\mapsto \mathscr{SB}$ and $K\mapsto
    \mathbb{HSP}(\Ad(K))$ for $K\in L(\underline{\mathsf{G}}_\mathrm{ref})$.
    Then $\iota$ is a lattice isomorphism.  Furthermore, a variety in $L(\mathscr{A}_\mathrm{ref})$
    is finitely axiomatized \up(finitely generated as a variety\up) if and only if it is the image
    under $\iota$ of either $\mathsf{\underline{S}}$ or a finitely axiomatized \up(finitely
    generated, respectively\up) uH class of reflexive graphs.
\end{thm}

We prove Theorems~\ref{thm:fund} and~\ref{thm:fundref} in the next two sections.

\section{Proof of Theorem~\ref{thm:fund}}
\subsection{Equations satisfied by adjacency
semigroups}\label{subsec:equational} The variety of semigroups
generated by the class of Rees matrix semigroups over trivial
groups is reasonably well understood: it is generated by a 5-element
semigroup usually denoted by $A_{2}$ (see \cite{leevol} for example).
(In context of this paper $A_{2}$ can be thought as the
semigroup reduct of the adjacency semigroup $\A(\Sg_{2})$ where
$\Sg_2$ is the 2-vertex graph from Example~\ref{eg:simplegraphgen}.)
This semigroup was shown to have a finite identity basis by Trahtman~\cite{trah},
who gave the following elegant description of the identities: an identity
$u\approx v$ (where $u$ and $v$ are semigroup words) holds in $A_{2}$
if and only if $u$ and $v$ start with the same letter, end with the same letter
and share the same set of two letter subwords.  Thus the equational
theory of this variety corresponds to pairs of words having the
same ``adjacency patterns'', in the sense that a two letter
subword $xy$ records the fact that $x$ occurs next to (and before)
$y$.  This adjacency pattern can also be visualized as a graph on
the set of letters, with an edge from $x$ to $y$ if $xy$ is a
subword, and two distinct markers indicating the first and last
letters respectively.

In this subsection we show that the equational theory of $\mathscr{A}$
has the same kind of property with respect to a natural unary semigroup
notion of adjacency.  The interpretation is that each letter has two
sides---left and right---and that the operation $'$ reverses these.
A subword $xy$ corresponds to the right side of $x$ matching the left
side of $y$, while $x'y$ or any subword $(x\ldots )'y$ corresponds to
the left side of $x$ matching the left side of $y$.
To make this more precise, we give an inductive definition.
Under this definition, each letter $x$ in a word will have two associated
vertices corresponding to the left and right side.  The graph will
have an initial vertex, a final vertex as well as a set of (directed)
edges corresponding to adjacencies.

Let $u$ be a unary semigroup word, and $X$ be the alphabet of letters
appearing in $u$.  We construct a graph $G[u]$ on the set
$$\{\ell_{x}\mid x\in X\}\cup\{r_{x}\mid x\in X\}$$
with two marked vertices. If $u$ is a single letter (say $x$),
then the edge set (or \emph{adjacency set}) of $G[u]$ is empty.
The \emph{initial vertex} of a single letter $x$ is $\ell_{x}$ and the
\emph{final} (or \emph{terminal}) \emph{vertex} is $r_{x}$.

If $u$ is not a single letter, then it is of the form $v'$ or
$vw$ for some unary semigroup words $v,w$.  We deal
with the two cases separately.  If $u$ is of the form $v'$, where $v$
has set of adjacencies $S$, initial vertex $p_{a}$ and final vertex $q_{b}$ (where
$\{p,q\}\subseteq \{\ell,r\}$ and $a,b$ are letters appearing in $v$), then
the set of adjacencies of $u$ is also $S$, but the initial vertex of $u$ is
equal to the final vertex $q_{b}$ of $v$ and the final vertex of $u$ is equal
to the initial vertex $p_{a}$ of $v$.

Now say that $u$ is of the form $vw$ for some
unary semigroup words $v,w$, with adjacency set $S_{v}$ and $S_{w}$
respectively and with initial vertices $p_{a_{v}}$, $p_{a_{w}}$
respectively and final vertices $q_{b_{v}}$ and $q_{b_{w}}$
respectively. Then the adjacency set of $G[u]$ is $S_{v}\cup
S_{w}\cup\{(q_{b_{v}},p_{a_{w}})\}$, the initial vertex is $p_{a_{v}}$
and the final vertex is $q_{b_{w}}$.  Note that the word $u$ may be
broken up into a product of two unary words in a number of different
ways, however it is reasonably clear that this gives rise to the same
adjacency set and initial and final vertices (this basically
corresponds to the associativity of multiplication).

For example the word $a'(baa')'$ decomposes as $a'\cdot (baa')'$,
and so has initial vertex equal to the initial vertex of $a'$,
which in turn is equal to the terminal vertex of $a$, which is
$r_{a}$.  Likewise, its terminal vertex should be the terminal
vertex of $(baa')'$, which is the initial vertex of $baa'$, which
is $\ell_{b}$. Continuing, we see that the edge set of the
corresponding graph has edges $\{(\ell_{a},\ell_{a}),
(r_{a},r_{a}), (r_{b},\ell_{a})\}$.  This graph is the first graph
depicted in Figure~\ref{fig:egs} (the initial and final vertices
are indicated by a sourceless and targetless arrow respectively).
The second is the graph of either of the words $a(bc)'$ or $(b(ac')')'$.
The fact that $G[a(bc)']=G[(b(ac')')']$ will be of particular importance
in constructing a basis for the identities of $\mathscr{A}$.

\setlength{\unitlength}{1.5pt}
\begin{figure}[ht]
    \begin{picture}(190,60)
            \put(10,0){$G[a'(baa')']$}
    \put(0,10)
    {
        \begin{picture}(40,50)

            \put(10,20){\circle*{4}}
            \put(-1,39){$\ell_{a}$}
            \put(10,40){\circle*{4}}
            \put(-1,10){$r_{a}$}

            \put(40,20){\circle*{4}}
            \put(28,39){$\ell_{b}$}
            \put(40,40){\circle*{4}}
            \put(30,10){$r_{b}$}

            \put(-1,20){\vector(1,0){9}}
            \put(40,40){\vector(1,0){10}}

            \put(40,20){\vector(-3,2){28}}

            \qbezier(10,40)(4,48)(10,48)
            \qbezier(10,40)(16,48)(10,48)

            \qbezier(10,20)(4,12)(10,12)
            \qbezier(10,20)(16,12)(10,12)
        \end{picture}
    }

            \put(130,0){$G[a(bc)']=G[(b(ac')')']$}
    \put(120,10)
    {
        \begin{picture}(70,80)

            \put(10,20){\circle*{4}}
            \put(-1,42){$\ell_{a}$}
            \put(10,40){\circle*{4}}
            \put(-1,10){$r_{a}$}

            \put(40,20){\circle*{4}}
            \put(28,42){$\ell_{b}$}
            \put(40,40){\circle*{4}}
            \put(30,10){$r_{b}$}

            \put(70,20){\circle*{4}}
            \put(58,42){$\ell_{c}$}
            \put(70,40){\circle*{4}}
            \put(70,10){$r_{c}$}

            \put(-1,20){\vector(1,0){9}}
            \put(40,40){\vector(1,0){10}}

            \put(40,20){\vector(3,2){29}}
            \put(69.1,15){\vector(1,3){1}}
            \qbezier(10,20)(10,5)(40,5)
            \qbezier(70,20)(70,5)(40,5)
        \end{picture}
    }
\end{picture}
\caption{Two examples of graphs of unary words}\label{fig:egs}
\end{figure}
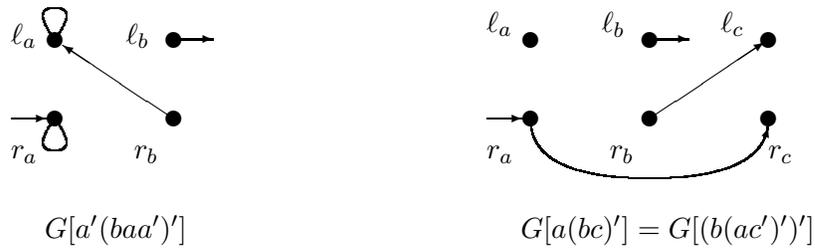

We can also construct a second kind of graph from a word $w$, in
which all loops are added to the graph of $G[w]$ (that is, it is the
reflexive closure of the edge set), we call this $G_\mathrm{ref}[w]$.
For example, it is easy to see that
$G_\mathrm{ref}[a'(baa')']=G_\mathrm{ref}[(ba)']$
(most of the work was done in the previous example).  Lastly, we
define the graph $G_\mathrm{symm}[w]$ corresponding to the symmetric
closure of the edge set of $G[w]$.

\begin{notn}\label{notation}
Let $u$ be a unary semigroup word and let $\theta$ be an assignment
of the letters of $u$ into nonzero elements of an adjacency semigroup
$\A(\Hg)$; say $\theta(x)=(i_{x},j_{x})$ for each letter $x$. Note that
$\theta(x')=(j_{x},i_{x})$, so we use the notation $i_{x'}:=j_{x}$
and $j_{x'}=i_{x}$.
\end{notn}

\begin{lem}\label{lem:lr}
Let $u$ and $\theta$ be as in Notation~\emph{\ref{notation}}.
If $\lambda_{a}$ is the initial vertex of $G[u]$ and $\rho_{b}$
is the terminal vertex \up(so $\lambda,\rho\in \{\ell,r\}$ and $a$
and $b$ are letters in $u$\up) and $\theta(u)\neq 0$, then
$\theta(u)=(i_{\bar{a}},j_{\bar{b}})$, where
$$\bar{a}=\begin{cases}
        a&\mbox{ if }\lambda=\ell\\
        a'&\mbox{ if }\lambda=r
    \end{cases}\quad\text{and}\quad
\bar{b}=\begin{cases}
        b&\mbox{ if }\rho=r\\
        b'&\mbox{ if }\rho=\ell.
    \end{cases}$$
\end{lem}

\begin{proof}
    This follows by an induction following the inductive definition
    of the graph of $u$.
\end{proof}

\begin{lem}\label{lem:graphhom}
Let $u$ and $\theta$ be as in Notation~\emph{\ref{notation}}.
Then $\theta(u)\neq 0$ if and only if the map defined by
$\ell_{x}\mapsto i_{x}$ and $r_{x}\mapsto j_{x}$ is a graph
homomorphism from $G[u]$ to $\Hg$.
\end{lem}
\begin{proof}
    Throughout the proof we use the notation of Lemma \ref{lem:lr}.

    (Necessity.) Say $\theta(u)\neq 0$, and let $(\rho_{b},\lambda_{a})$ be an edge in
    $G[u]$, where $\rho,\lambda\in\{\ell,r\}$ and $a$ and $b$ are letters in
    $u$).  We use Lemma \ref{lem:lr} to show that $(j_{\bar{b}},i_{\bar{a}})$ is
    an edge of $\Hg$. Note that in the case where no applications of $'$ are used
    (so we are dealing in the nonunary case), the  edge $(\rho_{b},\lambda_{a})$
    will necessarily be
    $(r_{b},\ell_{a})$; and we would want
    $(j_{b},i_{a})$ to be an edge of $\Hg$.

    Now, since $(\rho_{b},\lambda_{a})$ is an
    edge in $G[u]$, some subwords of $u$---say $u_{1}$ and $u_{2}$---have
    $u_{1}u_{2}$ a subword of $u$, and $\rho_{b}$ the terminal vertex of
    $G[u_{1}]$, and $\lambda_{a}$ the initial vertex of $G[u_{2}]$.
    Applying Lemma \ref{lem:lr} to both $G[u_{1}]$ and $G[u_{2}]$, we
    find that $\theta(u_{1})$ has right coordinate $j_{\bar{b}}$, and
    $\theta(u_{2})$ has left coordinate $i_{\bar{a}}$.  But $u_{1}u_{2}$ is a
    subword, so $\theta(u_{1})\theta(u_{2})\neq 0$, whence
    $(j_{\bar{b}},i_{\bar{a}})$ is an edge of $\Hg$, as required.

    (Sufficiency.)  This is easy.
\end{proof}
Lemma \ref{lem:graphhom} is easily
adapted to the graph $G_\mathrm{ref}(u)$ or $G_\mathrm{symm}(u)$, where the
graph $\Hg$ is assumed to be reflexive or symmetric, respectively.
\begin{pro}\label{pro:graphid}
    An identity $u\approx v$ holds in $\mathscr{A}$ if and only if $G[u]=G[v]$.  An
    identity holds in $\mathscr{A}_\mathrm{ref}$ if and only if $G_\mathrm{ref}[u]=G_\mathrm{ref}[v]$.
    An identity holds in $\mathscr{A}_\mathrm{symm}$ if and only if $G_\mathrm{    symm}[u]=G_\mathrm{symm}[v]$.
    \end{pro}
    \begin{proof}
    We prove only the first case; the other two cases are
    similar.

    First we show sufficiency.  Let us assume that $G[u]=G[v]$,
    and consider an assignment $\theta$ into an adjacency semigroup
    $\A(\Hg)$.  Now the vertex sets are the same, so $u$ and
    $v$ have the same alphabet.  So we may assume that $\theta$
    maps the alphabet to nonzero elements of $\A(\Hg)$.  By
    Lemma~\ref{lem:graphhom}, we have $\theta(u)\neq 0$ if
    and only if $\theta(v)\neq 0$.  By Lemma \ref{lem:lr}, we have
    $\theta(u)=\theta(v)$ whenever both sides are nonzero.
    Hence $\theta(u)=\theta(v)$ always.

    Now for necessity.  Say that $G[u]\neq G[v]$.  If the vertex
    sets are distinct, then $u\approx v$ fails on
    $\Ad(\underline{\bf 1})$, which is isomorphic to the unary
    semigroup formed by the integers $0$ and $1$ with the usual
    multiplication and the identity map as the unary operation.  Now
    say that $G[u]$ and $G[v]$ have the same vertices.  Without loss of
    generality, we may assume that either $G[v]$ contains an edge
    not in $G[u]$, or that the two graphs are identical but have
    different initial vertices.
    Let $\A_{u}:=\A(G[u])$ and consider the assignment into
    $\A_{u}$ that sends each variable $x$ to $(\ell_{x},r_{x})$.
    Observe that the value of $u$ is equal to $(\lambda_{a},\rho_{b})$
    where $\lambda_{a}$ is the initial vertex of $G[u]$ and $\rho_{b}$
    is the final vertex, while the value of $v$ is either $0$ (if there
    is an adjacency not in $G[v]$: we fail to get a graph homomorphism)
    or has different first coordinate (if $G[v]$ has a different initial vertex).
    So $u\approx v$ fails in $\mathscr{A}$.
    \end{proof}

\subsection{A normal form}\label{subsec:normalform}
Proposition \ref{pro:graphid} gives a reasonable solution to the
word problem in the $\mathscr{A}$-free algebras.  In this subsection
we go a bit further and show that every unary semigroup word is
equivalent in $\mathscr{A}$ to a unary semigroup word of a certain
form.  Because different forms may have the same adjacency graph,
this by itself does not constitute a different solution to the
word problem in $\mathscr{A}$-free algebras, however it is useful
in analyzing identities of $\mathscr{A}$.

Most of the work in this section revolves around the variety of
algebras of type $\langle 2,1\rangle$ defined by the three laws:
\[
\Psi=\{x''\approx x,\ x(yz)'\approx (y(xz')')',\ (xy)'z\approx
((x'z)'y)'\}.
\]
as interpreted within the variety of unary semigroups.   By examining
the adjacency graphs, it is easy to see that these identities are all
satisfied by $\mathscr{A}$ (see Fig.\,\ref{fig:egs} for one of
these).  In fact $\Psi$ defines a strictly larger variety than
$\mathscr{A}$ (it contains all groups for example), but they are close
enough for us to obtain useful information. For later reference we refer
to the second and third laws in $\Psi$ as the \emph{first associativity
across reversion} law (FAAR) and \emph{second associativity across reversion} law (SAAR),
respectively.  We let $\mathscr{B}$ denote the unary semigroup variety defined by
$\Psi$.

Surprisingly, the laws in $\Psi$ are sufficient to reduce every unary
semigroup word to one in which the nesting height of the unary $'$ is
at most $2$.  The proof of this is the main result of this subsection.

\begin{lem}\label{lem:bend}
     $\Psi$ implies $(a(bcd)'e)'\approx (b'e)'c(ad')'$, where $c$ is
     possibly empty.
\end{lem}
\begin{proof*}
    We prove the case where $c$ is non-empty only.  We have
    \[
     (a(bcd)'e)\fapprox
    ([bc(ad')']'e)'\sapprox
    ((b'e)'c(ad')')''\approx (b'e)'c(ad')'.\eqno{\Box}
    \]
\end{proof*}

Let $X:=\{x_{1},x_{2},\ldots\}$.  Let ${\bf F}(X)$
denote the free unary semigroup freely generated by $X$ and ${\bf
F}_{\Psi}(X)$ denote the $\mathscr{B}$-free algebra freely generated by $X$.
We let $\psi$ denote the fully invariant congruence on ${\bf F}(X)$ giving
${\bf F}_{\Psi}(X)={\bf F}(X)/\psi$.  We find a subset $N\subseteq {\bf
F}(X)$ with $X\subseteq N$ and show that multiplying two words from $N$
in ${\bf F}(X)$, or applying $'$ to a word in $N$ produces a word that is
$\psi$-equivalent to a word in $N$. In other words, $N$ forms a transversal of $\psi$;
equivalently, it shows that every word in ${\bf F}(X)$ is $\psi$-equivalent to a word in $N$.
In this way the members of $N$ are a kind of weak normal form for terms modulo $\Psi$
(we do not claim that distinct words in $N$ are not $\psi$-equivalent; for example,
Proposition \ref{pro:graphid} shows that $\mathscr{A}\models x(x'y)'\approx x(x'y')'$,
but the two words are distinct elements of $N$).

We let $N$ consist of all (nonempty) words of the form
\[
u_{1}(v_{1})'u_{2}(v_{2})'\ldots u_{n}(v_{n})'u_{n+1}
\]
for some some $n\in \omega$, where for $i\leq n$,
\begin{itemize}
    \item the $v_{i}$ are \emph{semigroup words} in the alphabet
    $$X\cup X'=\{x_{1},x_{1}',x_{2},x_{2}',\ldots\},$$ and all have
    length at least $2$ as semigroup words;
    \item the $u_{i}$ are possibly empty semigroup words in the
    alphabet $X\cup X'$ and if $n=0$, then $u_{1}$ is non-empty.
\end{itemize}
Notice that $X\subseteq N$ since the case $n=0$ corresponds to
semigroup words over $X\cup X'$.  For a member $s$ of $N$, we refer to
the number $n$ in this definition as the \emph{breadth} of $s$.

The following lemma is trivial.
\begin{lem}\label{lem:mult}
    If $s$ and $t$ be two words in $N$, then $s\cdot t$ is $\psi$-equivalent to a word in $N$.
\end{lem}
\begin{lem}\label{lem:dash}
    If $s$ is a word in $N$, then $s'$ is $\psi$-equivalent to a word in $N$.
\end{lem}
\begin{proof}
    We prove the lemma by induction on the breadth of $s$.
    If the breadth of $s$ is $0$ then $s'=(u_{1})'$ is either in $N$,
    or is of the form $x''$ for some variable $x$, in which case it reduces
    to $x\in N$ modulo $\Psi$.  Now say  that the result holds for breadth $k$
    members of $N$, and say  that the breadth of $s$ is $k+1$.  So $s$ can be
    written in the form $p(y_{1}\cdots y_{m})'u$ where $p$ is either empty
    or is a word from $N$ of breadth $k$,  $u=u_{k+2}$ is a possibly empty
    semigroup word in the alphabet $X\cup X'$ and $y_{1},\dots,y_{m}$ is a
    possibly repeating  sequence of variables from $X\cup X'$ with
    $v_{k+1}\equiv y_{1}\cdots y_{m}$ (so $m>1$).  Note that $p$ can be empty
    only if $k=0$.

    Let us write $w$ for $y_{2}\cdots y_{m-1}$ (if $m=2$, then $w$ is empty).
    If both $p$ and $u$ are empty, then $s'\in N$ already.  If neither $p$ nor $u$
    are empty, then by Lemma~\ref{lem:bend} $\Psi$ implies     $s'\approx (py_{m}')'w(y_{1}'u)'$.
    The breadth of $py_{m}'$ is $k$, so the induction hypothesis and Lemma~\ref{lem:mult}
    complete the proof.

    Now say that $p$ is  empty and $u$ is not.  We have
    $((y_{1}wy_{m})'u)'\sapprox (y_{1}'u)'wy_{m}$, and the latter word is contained
    in $N$ (modulo $x''\approx x$).

    Lastly, if $u$ is empty and $p$ is not, then we have $s'\equiv (p(wy_{m})')'\fapprox w(py_{m}')'$,
    and the  induction hypothesis applies to $(py_{m}')'$ since $py_{m}'$ is
    of breadth $k$. By Lemma~\ref{lem:mult}, $s'$ is $\psi$-equivalent to a member of $N$.
\end{proof}

As explained above, Lemmas \ref{lem:mult} and \ref{lem:dash} give us
the following result.
\begin{pro}\label{pro:nf}
    Every unary semigroup word reduces modulo $\Psi$ to a word in $N$.
\end{pro}
A algorithm for
making such a reduction is to iterate the method of
proof of Lemma \ref{lem:mult} and \ref{lem:dash}, however we will not need this here.

\subsection{Subvarieties of $\mathscr{A}$ and sub-uH classes of
$\underline{\mathsf{G}}$}\label{subsec:fund}
In this subsection we complete the proof of Theorem~\ref{thm:fund}. Recall
that the theorem claims that, for any nonempty class $K$ of graphs, any
graph $\Gg$ belongs to the uH-class generated by $K$ if and only
if the adjacency semigroup $\A(\Gg)$ belongs to the variety generated by
the adjacency  semigroups $\A(\Hg)$ with $\Hg\in K$. For the ``only if''
statement we use a direct argument.  For the ``if'' statement, we use
a syntactic argument, translating uH sentences of $K$ into identities
of $\A(K)$.

\begin{lem}\label{lem:uHtoV}
    If $\Gg\in\mathbb{ISP^{+}P}_\mathrm{u}(K)$, then $\Ad(\Gg)\in\mathbb{HSP}(\Ad(K))$.
\end{lem}
\begin{proof}
    First consider a nonempty family $L=\{\Hg_{i}\mid
    i\in I\}$ of graphs from $K$ and an ultraproduct
    $\Hg:=\prod_{U}L$ (for some ultrafilter $U$ on $2^{I}$).
    It is easy to see that the ultraproduct of the family
    $\{\A(\Hg_{i})\mid i\in I\}$ over the same ultrafilter
    $U$ is isomorphic to $\A(\Hg)$ (we leave this elementary
    proof to the reader).  Hence, we have $\mathbb{I}(\A\mathbb{P}_\mathrm{u}(K))=\mathbb{IP}_\mathrm{u}(\A(K))$.
    Now we have $\Gg\in\mathbb{ISP^{+}}(\mathbb{P}_\mathrm{u}(K))$.  So it
    will suffice to prove that $\A(\Gg)\in\mathbb{HSP}(\A(\mathbb{P}_\mathrm{u}(K)))$,
    since $\mathbb{HSP}(\A(\mathbb{P}_\mathrm{u}(K)))= \mathbb{HSPP}_\mathrm{u}(\A(K))=\mathbb{HSP}(\A(K))$.  We let
    $P$ denote $\mathbb{P}_\mathrm{u}(K)$.

    Now $\Gg$ is isomorphic to an induced subgraph of the direct product
    $\prod_{i\in I}\Hg_{i}$ with $\Hg_{i}\in P$.  It does no
    harm to assume that this embedding is the inclusion map.  Let
    $\pi_{i}:\Gg\to \Hg_{i}$ denote the projection.
    Evidently the following properties hold:
    \begin{itemize}
    \item[(i)] if $u$ and $v$ are distinct vertices of $\Gg$
    then there is $i\in I$ such that $\pi_{i}(u)\neq \pi_{i}(v)$;
    \item[(ii)] if $(u,v)$ is not an edge of $\Gg$ then there
    is $i\in I$ with $(\pi_{i}(u),\pi_{i}(v))$ not an edge of $\Hg_{i}$.
    \end{itemize}
    We aim to show that $\A(\Gg)$ is a quotient of a subalgebra of $\prod_{i\in
    I}\A(\Hg_{i})$.  We define a map $\alpha:\A(\Gg)\to
    \prod_{i\in I}\A(\Hg_{i})$ by letting $\alpha(0)$ be the constant $0$ and
    $\alpha(u,v)$ be the map $i\mapsto (\pi_{i}(u),\pi_{i}(v))$.  The
    map $\alpha$ is unlikely to be a homomorphism.  Let $B$ be
    the subalgebra of $\prod_{i\in
    I}\A(\Hg_{i})$ generated by the image of
    $\A(\Gg)$, and let $J$ be the ideal of $B$ consisting of all
    elements with a $0$ coordinate.

    \textbf{Claim 1.} Say $(u_{1},v_{1})$ and $(u_{2},v_{2})$ are
    (nonzero) elements of $\A(\Gg)$.  If $v_{1}\sim u_{2}$ then
    $\alpha(u_{1},v_{1})\alpha(u_{2},v_{2})=\alpha((u_{1},v_{1})(u_{2},v_{2}))$.
    \begin{proof}
    Now
    $\alpha((u_{1},v_{1})(u_{2},v_{2}))[i]=\alpha(u_{1},v_{2})[i]
    =(\pi_{i}(u_{1}),\pi_{i}(v_{2}))$, because $v_{1}\sim u_{2}$
    in $\Gg$.  Also, for every $i\in I$ we
    have $\pi_{i}(v_{1})\sim\pi_{i}(u_{2})$, so that $\alpha(u_{1},v_{1})[i]\alpha(u_{2},v_{2})[i]=
    (\pi_{i}(u_{1}),\pi_{i}(v_{1}))(\pi_{i}(u_{2}),\pi_{i}(v_{2}))=(\pi_{i}(u_{1}),
    \pi_{i}(v_{2}))$ as required.
    \end{proof}

    \textbf{Claim 2.}  Say $(u_{1},v_{1})$ and $(u_{2},v_{2})$ are
    nonzero elements of $\A(\Gg)$.  If $v_{1}\not\sim u_{2}$ then
    $\alpha(u_{1},v_{1})\alpha(u_{2},v_{2})\in I$.
    \begin{proof}
    By the definition of $v_{1}\not\sim u_{2}$ there is $i\in I$ with
    $\pi_{i}:\Gg\to \Hg_i$ with $\pi_{i}(v_{1})\not\sim\pi_{i}(u_{2})$.
    Then
    $(u_{1},v_{1})[i](u_{2},v_{2})[i]=0$.
    \end{proof}

    Claims 1 and 2 show that $\alpha$ is a semigroup homomorphism
    from $\A(\Gg)$ onto $B/I$ (at least, if we adjust the co-domain of $\alpha$ to
    be $B/I$ and identify the constant $0$ with $I$).  Now we show that this map is injective.
    Say $(u_{1},v_{1})\neq (u_{2},v_{2})$ in $\A({\Gg})$.  Without
    loss of generality, we may assume that $u_{1}\neq u_{2}$.  So
    there is a coordinate $i$ with $\pi_{i}(u_{1})\neq
    \pi_{i}(u_{2})$.  Then $\alpha(u_{1},v_{1})$ differs from
    $\alpha(u_{2},v_{2})$ on the $i$-coordinate.  So we have a
    semigroup isomorphism from $\A(\Gg)$ to $B/I$.  Lastly, we observe
    that $\alpha$ trivially preserves the unary operation, so we have
    an isomorphism of unary semigroups as well.  This completes
    the proof of Lemma \ref{lem:uHtoV}.
\end{proof}

To prove the other half of Theorem \ref{thm:fund} we take a syntactic
approach by translating uH sentences into unary semigroup identities.
To apply our technique, we first need to reduce arbitrary uH sentences
to logically equivalent ones of a special form.

Our goal is to show that if $\Gg\notin \mathbb{ISP^{+}P}_\mathrm{u}(K)$ then
$\A(\Gg)\notin\mathbb{HSP}(\A(K))$.  We first consider some degenerate
cases.

If $K=\{\underline{\mathbf{0}}\}$, then $\A(K)$ is the class
consisting of the one element unary semigroup and $\mathbb{HSP}(\A(K))\models x\approx y$.
The statement $\Gg\notin \mathbb{ISP^{+}P}_\mathrm{u}(K)$ simply means that $|G|\geq 1$ and so $\A(\Gg)\not\models
x\approx y$.  So $\A(\Gg)\notin \mathbb{HSP}(K)$.

Now we say that $K$ contains a nonempty graph.  We can
then further assume that the empty graph is not in $K$.  If
$\Gg$ is the 1-vertex looped graph $\underline{\bf 1}$, then
the statement $\Gg\notin \mathbb{ISP^{+}P}_\mathrm{u}(K)$ simply
means that $K$ consists of antireflexive graphs.  In this case,
$\A(K)\models xx'\approx 0$, while $\A(\Gg)\not\models xx'\approx 0$.
So again, $\A(\Gg)\notin \mathbb{HSP}(\A(K))$.

So now it remains to consider the case where $\Gg$ is not the
1-vertex looped graph and $K$ does not contain the empty
graph.  Lemma \ref{lem:malcev} shows that there is some uH
sentence $\Gamma$ holding in each member of $K$, but failing
on $\Gg$. We now show that $\Gamma$ can be chosen to be a
quasi-identity.

If $\Gamma$ is a uH sentence of
the second kind, say $\bigvee_{1\leq i\leq n}\neg \Phi_{i}$, then
choose some atomic formula $\Xi$ in variables not appearing in any of
the $\Phi_{i}$ and that fails on $\Gg$ under some assignment: for
example, if $|G|\geq 2$, then a formula of the form $x\approx y$
suffices, while if $|G|=1$, then $\Gg$ is the one element
edgeless graph and $x\sim y$ suffices.  Now
replace $\Gamma$ with the quasi-identity $\bigand_{1\leq i\leq
n}\Phi_{i}\rightarrow \Xi$.  We need to show that
this new quasi-identity holds in $K$ and fails in $\Gg$.
It certainly holds in $K$, since it is logically
equivalent to $\bigvee\neg \Phi_{i}\vee \Xi$, while $\bigvee\neg
\Phi_{i}$ is constantly true.  On the other hand,
since  $\bigvee_{1\leq i\leq n}\neg \Phi_{i}$ does not hold on $\Gg$,
there is an assignment $\theta$ making $\bigand_{1\leq i\leq
n}\Phi_{i}$ true, and this assignment can be extended to the
variables of $\Xi$ in such a way that $\Xi$ is false under $\theta$.
In other words, we have a failing assignment for the new
quasi-identity on $\Gg$.

Next we need to show that the quasi-identity $\Gamma$ can be chosen
to have a particular form.
Let us call a quasi-identity \emph{reduced}
if the equality symbol $\approx$ does not appear
in the premise.  One may associate any quasi-identity with a reduced
quasi-identity in the obvious way: if $x\approx y$ appears in the
premise of the original, then we may replace all occurrences of $y$
in the quasi-identity with $x$ (including in the conclusion) and remove $x\approx y$ from the
premise.  If a quasi-identity fails on a graph $\Hg$ under some
assignment $\theta$, then the corresponding reduced
quasi-identity also fails under $\theta$.  Conversely, if the
reduced quasi-identity fails on $\Hg$ under some assignment
$\theta$, then we may extend $\theta$ to a failing assignment of the
original quasi-identity.  This means that we may choose $\Gamma$ to
be a reduced quasi-identity.

Let $\Phi:=\&_{1\leq i\leq n} u_{i}\sim v_{i}$ be a conjunction of
adjacencies, where the
\[
u_{1},\ldots,u_{n},v_{1},\ldots,v_{n}\in\{a_{1},\ldots,a_{m}\}
\]
are not necessarily distinct variables.  For each adjacency $u_{i}\sim v_{i}$ in $\Phi$,
let $w_{i}$ denote the word
$(u_{i}v_{i})'s_{i}(u_{i}v_{i}')'s_{i}(u_{i}'v_{i})'s_{i}(u_{i}'v_{i}')'$,
where $s_{i}$ is a new variable.  Now let
$\sigma:\{1,\ldots,m\}\to\{1,\ldots,n\}$ be some finite sequence of numbers from
$\{1,\ldots,n\}$ with the property that for each pair $i,j\in
\{1,\ldots,n\}$ there is $k< m$ with $\sigma_{k}=i$ and
$\sigma_{k+1}=j$ and such that $\sigma(1)=\sigma(m)=1$.  Define a word
$D_{\Phi}$ (depending on $\sigma$) as follows:
\[
\left(\prod_{1\leq
i<m}w_{\sigma(i)}t_{\sigma(i),\sigma(i+1)}\right)w_{\sigma(m)},
\]
where the $t_{i,j}$ are new variables.
As an example, consider the conjunction $\Phi:=x\sim y\And y\sim z$,
where $n=2$, $u_{1}=x$, $v_{1}=u_{2}=y$ and $v_{2}=z$.  Using the
sequence $\sigma=1,2,2,1,1$ we get $D_{\Phi}$ equal to the following
expression:
\begin{multline*}
(xy)'s_{1}(xy')'s_{1}(x'y)'s_{1}(x'y')'t_{1,2}
(yz)'s_{2}(yz')'s_{2}(y'z)'s_{2}(y'z')'t_{2,2}\\
(yz)'s_{2}(yz')'s_{2}(y'z)'s_{2}(y'z')'t_{2,1}\\
(xy)'s_{1}(xy')'s_{1}(x'y)'s_{1}(xy)'t_{1,1}
(xy)'s_{1}(xy')'s_{1}(x'y)'s_{1}(x'y')'.
\end{multline*}
\begin{lem}\label{lem:homom}
    Let $\Hg$ be a graph, $\Phi$ be a conjunction of adjacencies in
    variables $a_{1},\ldots,a_{m}$ and $\theta$ be
    an assignment of the variables
    of $D_{\Phi}$ into $\A(\Hg)$\up, with
    $\theta(a_{i})=(\ell_{i},r_{i})$ say.
    Let $\gamma$ be any member of
    $\{L,R\}^{m}$.  If $\theta(D_{\Phi})\neq 0$ then the map
    $\phi_{\gamma}$ from $a_{1},\ldots,a_{m}$ into the vertices
    $V_\Hg$ of $\Hg$ defined by
    \[
    \phi_{\gamma}(a_{i})=
    \begin{cases}
    \ell_{i}&\mbox{  if }\gamma(i)=L;\\
    r_{i}&\mbox{  if }\gamma(i)=R
    \end{cases}
    \]
    satisfies $\Phi$.
    \end{lem}
    \begin{proof}
    Let $a_{i}\sim a_{j}$ be one of the adjacencies in $\Phi$.
    So all of $a_{i}a_{j}$, $a_{i}'a_{j}$, $a_{i}a_{j}'$
    and $a_{i}'a_{j}'$ appear in $D_{\Phi}$ and hence are given
    nonzero values by $\theta$. We have
    $\ell_{i}\sim\ell_{j}$, $\ell_{i}\sim r_{j}$, $r_{i}\sim \ell_{j}$,
    $r_{i}\sim r_{j}$ in $\Hg$. So regardless of the choice of
    $\gamma$ we have $\phi_{\gamma}(a_i)\sim \phi_{\gamma}(a_j)$ in $\Hg$.
    \end{proof}
    \begin{lem}\label{lem:plus}
    Let $\Phi=\bigand_{1\leq i\leq n}u_{i}\sim v_{i}$ be a
    nonempty conjunction in the variables $a_{1},\ldots,a_{m}$
    and let $\theta$ be an assignment of these variables
    into a graph $\Hg$ such that $\Hg\models \theta(\Phi)$.
    Define an assignment $\theta^{+}$ of the
    variables of $D_{\Phi}$ into $\A(\Hg)$ by $a_{i}\mapsto
    (\theta(a_{i}),\theta(a_{i}))$,
    $\theta^{+}(t_{i,j}):=(\theta(v_{i}),\theta(u_{j}))$ and
    $\theta^{+}(s_{i})=\theta^{+}(t_{i,i})$.  We have
    $\theta^{+}(D_{\Phi})=(\theta(v_{1}),\theta(u_{1}))$.
    \end{lem}
    \begin{proof}  This is a routine calculation.
    For each adjacency
    $u_{i}\sim v_{i}$ in $\Phi$ (here $\{u_{i},v_{i}\}\subseteq
    \{a_{1},\ldots,a_{m}\}$) we have
    $$\theta^{+}((u_{i}v_{i})'),\ \theta^{+}((u_{i}v_{i}')'),\
    \theta^{+}((u_{i}'v_{i})'),\ \theta^{+}((u_{i}'v_{i}')')$$ all
    taking the same nonzero value $(\theta(v_{i}),\theta(u_{i}))$.
    Then we also have
    $\theta^{+}(w_{i})=(\theta(v_{i}),\theta(u_{i}))$ which shows
    that
    \begin{multline*}
    \theta^{+}(D_{\Phi})\\=[\theta(u_{1}),\theta(v_{1})]\ldots
    [\theta(v_{i}),\theta(u_{i})]\ \theta^{+}(t_{i,j})\ [\theta(v_{j}),\theta(u_{j})]\ldots
    [\theta(u_{1}),\theta(v_{1})]\\
    =[\theta(u_{1}),\theta(v_{1})]\ldots
    [\theta(v_{i}),\theta(u_{i})][\theta(v_{i}),\theta(u_{j})][\theta(v_{j}),\theta(u_{j})]\ldots
    [\theta(u_{1}),\theta(v_{1})]\\
    =[\theta(u_{1}),\theta(v_{1})]
    \end{multline*}
    (where the square brackets are used for clarity only).
    \end{proof}
    \begin{lem}\label{lem:qiequals}
    Let $\Hg$ be a nonempty graph and  $\Phi\rightarrow u\approx v$ be a
    reduced   quasi-identity where $\Phi$ is nonempty and one of $u$ or $v$
    does not appear in $\Phi$ (say it is $u$).  We have $\A(\Hg)\models uD_{\Phi}\approx
    u'D_{\Phi}$ if and only if $\Hg\models \Phi\rightarrow u\approx v$.
    \end{lem}
    \begin{proof}
    First assume that $\Hg\models \Phi\rightarrow u\approx
    v$, where $u$ does not appear in $\Phi$.
    Both sides of the identity contain the subword $D_{\Phi}$ and so we may
    consider an assignment $\theta$ sending $D_{\Phi}$ to a
    nonzero value (if there are none, then we are done).  By
    Lemma \ref{lem:homom}, we have an interpretation of $\Phi$ in
    $\Hg$.  But then, we can choose any value for $\theta(u)$ and
    find that it is the same value as $\theta(v)$.  In other
    words,  $\Hg$ has only one vertex.  Also, since $D_{\Phi}$ takes
    a nonzero value on $\A(\Hg)$, we find that the semigroup
    reduct of $\A(\Hg)$ is not a null semigroup (that is, a semigroup
    in which all products are equal to 0).  Hence $\A(\Hg)$ is isomorphic
    to the unary  semigroup formed by the integers $0$ and $1$ with the usual
    multiplication and the identity map as the unary operation. In this
    case we have $u\approx u'$ satisfied and the identity holds.

    Now say that $\Hg\not\models \Phi\rightarrow u\approx v$, and
    let $\theta$ be a failing assignment.  As
    $\theta(u)\neq \theta(v)$ we can find a
    vertex $a$ of $\Hg$ such that $a\neq \theta(u_{1})$.
    Extend the assignment $\theta^{+}$ of Lemma \ref{lem:plus} by
    $u\mapsto(a,\theta(u_{1}))$.  Evidently,
    $\theta^{+}(uD_{\Phi})=(a,\theta(u_{1}))(\theta(v_{1}),(\theta(u_{1}))=(a,\theta(u_{1}))$, but
    $\theta^{+}(u'D_{\Phi})$ is either equal to $0$, or is non
    zero but has left coordinate different to $\theta^{+}(uD_{\Phi})$.
    \end{proof}

    \begin{lem}\label{lem:qiequals2}
    Let $\Hg$ be a nonempty graph and  $\bigand_{1\leq i\leq
    n}u_{i}\sim v_{i}\rightarrow u\approx v$ be a
    reduced
    quasi-identity where $\Phi$ is nonempty and both $u$ and $v$
    appear in $\Phi$.  If $u=u_{i}$ for some $i$ then we have $\A(\Hg)\models
    u_{i}t_{i,1}D_{\Phi}\approx
    vt_{i,1}D_{\Phi}$ if and only if $\Hg\models \Phi\rightarrow
    u\approx v$.  If $u=v_{i}$ for some $i$ then  we have $\A(\Hg)\models
    D_{\Phi}t_{1,i}u_{i}\approx
    D_{\Phi}t_{1,i}v$ if and only if $\Hg\models \Phi\rightarrow
    u\approx v$
    \end{lem}
    \begin{proof}
    First assume that $\Hg\models \Phi\rightarrow u\approx
    v$.  We consider only the case that $u=u_{i}$; the other case
    follows by symmetry.  As before, we can consider the case
    where there is an assignment $\theta$ into $\Hg$
    satisfying $\Phi$.  So we have $\theta(u_{i})=\theta(v)$ and
    hence $\theta^{+}(u_{i})=\theta^{+}(v)$, in
    which case both sides of the identity take the same value.

    Now say that $\Hg\not\models \Phi\rightarrow u\approx
    v$ and let  $\theta$ be a failing assignment.  Now, the left
    side of the identity contains the same adjacencies as $D_{\Phi}$
    and so takes a nonzero value in $\A(\Hg)$ under
    the assignment $\theta^{+}$; moreover the left coordinate is
    $\theta(u_{i})$.  However the right hand side
    either takes the value $0$ (if $\theta(v)\not\sim
    \theta(v_{i})$) or has left coordinate equal to
    $\theta(v)\neq\theta(u_{i})$.  In either case, the identity
    fails.
    \end{proof}

    Now we come to reduced quasi-identities in which the conclusion
    is an adjacency $u\sim v$.    We consider $9$ cases according to
    whether or not $u$ and $v$ appear in $\Phi$, and if so, whether
    they appear  as the ``source'' or ``target'' of an adjacency.
    The nine identities $\tau_{1},\ldots,\tau_{9}$ are defined in the
    following table.  In
    this table, the first row corresponds to the situation
    where neither $u$ nor $v$ appear, while the second corresponds
    to the situation where $u$ does not appear, but $v$ does appear
    as some $u_{j}$ (in other words, as a ``source''), and so on.  If
    one of $u$ or $v$ appears as both a source and a target, then
    there will be choices as to which identity we can choose.
    The variables $z$ and $w$ are
    new variables not appearing in $D_{\Phi}$.
    \begin{center}
    \begin{tabular}{c|c|c}
        $k$&$u\sim v$&$\tau_{k}$\\
        \hline
        $1$.\rule{0cm}{1em} & $z\sim w$&
        $wD_{\Phi}z\approx
        (wD_{\Phi}z)^{2}$\\
        $2$.\rule{0cm}{1em}&
        $z\sim u_{j}$&
        $u_{j}t_{j,1}D_{\Phi}z\approx
        (u_{j}t_{j,1}D_{\Phi}z)^{2}$\\
        $3$.\rule{0cm}{1em}&
        $z\sim v_{j}$&
        $(u_{j}v_{j})'t_{j,1}D_{\Phi}z\approx
        ((u_{j}v_{j})'t_{j,1}D_{\Phi}z)^{2}$\\
        $4$.\rule{0cm}{1em}&
        $u_{i}\sim w$&
        $wD_{\Phi}t_{1,i}(u_{i}v_{i})'\approx
        (wD_{\Phi}t_{1,i}(u_{i}v_{i})')^{2}$\\
        $5$.\rule{0cm}{1em}&
        $v_{i}\sim w$&
        $wD_{\Phi}t_{1,i}v_{i}\approx
        (wD_{\Phi}t_{1,i}v_{i})^{2}$\\
        $6$.\rule{0cm}{1em}&
        $u_{i}\sim
        u_{j}$&$u_{j}t_{j,1}D_{\Phi}t_{1,i}(u_{i}v_{i})'\approx
        (u_{j}t_{j,1}D_{\Phi}t_{1,i}(u_{i}v_{i})')^{2}$\\
        $7$.\rule{0cm}{1em}&
        $u_{i}\sim
        v_{j}$&$(u_{j}v_{j})'t_{j,1}D_{\Phi}t_{1,i}(u_{i}v_{i})'\approx
        ((u_{j}v_{j})'t_{j,1}D_{\Phi}t_{1,i}(u_{i}v_{i})')^{2}$\\
        $8$.\rule{0cm}{1em}&
        $v_{i}\sim
        u_{j}$&$u_{j}t_{j,1}D_{\Phi}t_{1,i}v_{i}\approx
        (u_{j}t_{j,1}D_{\Phi}t_{1,i}v_{i})^{2}$\\
        $9$.\rule{0cm}{1em}&
        $v_{i}\sim v_{j}$&$(u_{j}v_{j})'t_{j,1}D_{\Phi}t_{1,i}v_{i}\approx
        ((u_{j}v_{j})'t_{j,1}D_{\Phi}t_{1,i}v_{i})^{2}$
    \end{tabular}
    \end{center}

    \begin{lem}\label{lem:qisim}
    Let $\Hg$ be a graph and  $\Phi\rightarrow u\sim v$ be a
    quasi-identity where $\Phi$ is nonempty.  Consider the
    corresponding identity $\tau_{k}$.  We have $\A(\Hg)\models
    \tau_{k}$ if and only if $\Hg\models
    \Phi\rightarrow u\approx v$.
    \end{lem}
    \begin{proof}
    We prove the case
    of $\tau_{4}$ and leave the remaining (very similar) cases to the reader.
    First assume that $\Hg\models \Phi\rightarrow u_{i}\sim
    v$.  Consider some assignment $\theta$  into $\A(\Hg)$
    that gives $D_{\Phi}$ a nonzero value.  As $w$ appears on
    both sides, we may further assume that $\theta(w)$ is nonzero.  Observe that the
    graph of the right hand side of the identity is identical to
    that of the left side except for the addition of a single
    edge from $\ell_{u_{i}}$ to $\ell_{w}$.  Also, the initial and final vertices are
    the same.  So to show that the two sides are equal, it
    suffices to show that $\theta(u_{i}')\theta(w)$ is non
    zero.

    Choose any map $\gamma$ from the variables of $\tau_{4}$
    to $\{L,R\}$ with $\gamma(u_{i})=R$.  By
    Lemma \ref{lem:homom} we have $\Hg\models\phi_{\gamma}(\Phi)$.  Using $\Phi\rightarrow
    u_{i}\sim v$ it follows that for any
    vertex $w$ we have $\phi_{\gamma}(u_{i})\sim w$.  In other
    words, $\theta(u_{i}')\theta(w)$ is nonzero as required.

    Now say that $\Phi\rightarrow u_{i}\sim v$ fails on $\Hg$ under
    some assignment $\theta$.  Extend $\theta^{+}$ to $w$ by
    $w\mapsto (\theta(v),\theta(u_{1}))$.  Under this assignment
    the left hand side of $\tau_{4}$ takes the value
    $(\theta(v),\theta(u_{i}))$, while the right hand side equals
    $0$.
    \end{proof}

    Lastly we need to consider the case where $\Gamma$ has empty
    premise, that is, where $\Gamma$ is a universally quantified
    atomic formula
    $\tau$.  In the language of
    graphs, there are essentially four different possibilities
    for $\tau$ (up
    to a permutation of letter names):
    $x\sim y$, $x\sim x$, $x\approx y$ and  $x\approx x$.  The
    last of these is a tautology.  The first three are
    nontautological and correspond to the uH-classes of complete
    looped graphs, reflexive graphs, and the one element graphs.
    For $\Phi$ one of the three atomic formulas, we let
    $\tau_{\Phi}$ denote the identities $xx\approx x$,
    $xx'x\approx x$, and $x'\approx x$, respectively.
    \begin{lem}\label{lem:atomic}
        Let $\Hg$ be a graph and $\Phi$ be one of the three
        nontautological atomic formulas in the language of
        graphs.  We have $\Hg\models \Phi$ if and only if
        $\A({\Hg})\models \tau_{\Phi}$.
    \end{lem}
    \begin{proof}
        If $\Phi$ is $x\sim y$, then it is easy to see that
        $\Hg\models \Phi$ if and only if the underlying
        semigroup of $\A({\Hg})$ satisfies $xx\approx x$.  The
        case of $\Phi=x\sim x$ has been discussed already in
        Section \ref{sec:main}.  The case of $x\approx y$ corresponds
        to the 1-vertex graphs, which is clearly equivalent to the property
        that $\A({\Hg})\models x'\approx x$.
    \end{proof}

    Now we can complete the proof of Theorem \ref{thm:fund}.  We
    have a reduced quasi-identity $\Gamma$ satisfied by
    $K$ and failing on $\Gg$.  By the appropriate
    choice out of Lemmas \ref{lem:qiequals}, \ref{lem:qiequals2},
    \ref{lem:qisim} or \ref{lem:atomic} we can construct an identity
    $\tau$ such  that $\A(K)\models \tau$
    and $\A({\Gg})\not\models \tau$.  Hence $\A({\Gg})\notin\mathbb{HSP}(\A(K))$.\hfill$\Box$

\section{Proof of Theorem \ref{thm:fundref}}
\label{sec:regular}

In contrast to the proof of Theorem \ref{thm:fund}, this section requires
some basic notions and facts from semigroup theory such as Green's relations $\mathrsfs{J}$,
$\mathrsfs{L}$, $\mathrsfs{R}$, $\mathrsfs{H}$ and their manifestation on Rees matrix
semigroups. For details, refer to the early chapters of any general semigroup theory
text; Howie \cite{how} for example.

The first step to proving Theorem \ref{thm:fundref} is the following.

\begin{lem}\label{lem:reflexive}
    Let $\mathscr{V}$ be a variety of unary semigroups satisfying
\begin{equation}
\label{eq:simple}
    xx'x\approx x,\ x''\approx x,\ (x'x)'\approx x'x,\ (xy)'\approx
    y'(x'xyy')'x'.
\end{equation}
    If $A\in\mathscr{V}$ as a semigroup is a
    completely $0$-simple semigroup with trivial subgroups,
    then $A$ is of the form $\A(\Hg)$ for some reflexive graph $\Hg$.
\end{lem}
\begin{proof}
    Since $A$ is a completely 0-simple semigroup with trivial subgroups,
    the Green relation $\mathrsfs{H}$ is trivial.  Now every $a\neq 0$ in
    $A$ is $\mathrsfs{L}$-related to $a'a$ (since $a(a'a)=a$) and
    $\mathrsfs{R}$-related to $aa'$.  Also, these elements are fixed by
    $'$ by identity $(x'x)'\approx x'x$ (and $x''\approx x$).  Next
    we observe that $a\Lgr b$ if and only if  $a'\Rgr b'$.  For this
    we can use identity $(xy)'\approx y'(x'xyy')'x'$: if $a\Lgr b$ then $xa=b$ for
    some $x$, so $b'=(xa)'=a'z$, for $z=(x'xaa')x'$.  So $b'\Rgr a'$.  The
    other case follows by symmetry (or using $(x')'\approx x$).

    This implies that each $\mathrsfs{L}$-class and each
    $\mathrsfs{R}$-class contain precisely one fixed point of~$'$ (if
    $a'=a$, $b'=b$ and $a\Lgr b$, then $a=a'\Rgr b'=b$, so
    $a\Hgr b$).  Represent  $A$ as a Rees matrix semigroup
    (with matrix $P$) in which fixed points of $'$
    correspond to diagonal elements (as $xx'$ is an idempotent,
    $P$ will have $1$ down the diagonal).  It is easily seen this is
    $\A(\Hg)$ for the graph $\Hg$ with $P$ as adjacency matrix.  This
    graph is reflexive since the identity $xx'x\approx x$ holds.
\end{proof}

In the case where $\Hg$ is a universal relation, the set
$\A(\Hg)\backslash\{0\}$ is a subuniverse, and the corresponding
subalgebra of $\A(\Hg)$ is a square band.

\begin{lem}\label{lem:reflexive2}
    Let $\mathscr{V}$ be a variety of unary semigroups satisfying
    the identities \eqref{eq:simple}.  If $A\in\mathscr{V}$
    as a semigroup is a  completely simple semigroup with trivial
    subgroups, then $A$ is a square band.
\end{lem}
\begin{proof}
    The proof is basically the same as for Lemma \ref{lem:reflexive}.
\end{proof}

In order to get a small basis for the identities of $\mathscr{A}_\mathrm{ref}$ the
following lemma is useful.
    \begin{lem}\label{lem:easycons}
    The following laws are consequences of
    $$\Psi_{1}:=\{x\approx xx'x,(x'x)'\approx x'x,x''\approx x,x(yz)'\approx (y(xz')')',
    (xy)'z\approx((x'z)'y)'\}{:}$$
    \begin{itemize}
     \item $(xy)'\approx y'(xyy')'\approx (x'xy)'x'\approx y'(x'xyy')'x'$\up;
    \item $(xyz)'\approx (yz)'y(xy)'$.
    \end{itemize}
    \end{lem}
    \begin{proof*}
    For the first item we have $\Psi_{1}$ implies $(xy)'\approx ((xx')'xy)'\sapprox (x'xy)'x'$.
    The other two cases of this item are very similar.

    For the second item, first note that using item 1 and $\Psi_{1}$, we have $(xyy')'\approx
    y(xyy'y)'\approx y(xy)'$.
    Using this we obtain
    $$(xyz)'{\approx}(xy(y'y)'z)'\fapprox
    ((y'(xyy')')'z)'\sapprox (yz)'(xyy')'{\approx}(yz)'y(xy)'.\eqno{\Box}$$
     \end{proof*}
    The second item of Lemma \ref{lem:easycons} enables a refinement of
    Proposition \ref{pro:nf}.
    \begin{cor}\label{cor:refnf}
    The identities $\Psi_1$ reduce every unary semigroup word
    to a member of $N$ in which each
    subword of the form $(v)'$ has the property that $v$ is
    a semigroup word of length $1$ or $2$ \up(over the alphabet $X\cup X'$\up{).}
    \end{cor}

    We now let $\Sigma_\mathrm{ref}$ denote the following set of unary
    semigroup identities:
    \begin{align}
    & x''\approx x,\ x(yz)'\approx
    (y(xz')')',\
    (xy)'z\approx ((x'z)'y)'\tag{$\Psi$}\\
    &xx'x\approx x,\label{eq:reg}\\
    &(xx')'\approx xx',\label{eq:fix}\\
    &x^{3}\approx x^{2},\label{eq:aperiodic}\\
    &xyxzx\approx xzxyxzx\approx xzxyx,\label{eq:m}\\
    &x'yxzx\approx (xzx)'yxzx,\label{eq:r}\\
    &xyxzx'\approx xyxz(xyx)'.\label{eq:l}
    \end{align}
    Proposition
    \ref{pro:graphid} easily shows that all but identity
    \eqref{eq:reg} hold in $\mathscr{A}$, while  \eqref{eq:reg}
    obviously holds in the subvariety $\mathscr{A}_\mathrm{ref}$.
    Hence, to prove that $\Sigma_\mathrm{ref}$ is a basis for $\mathscr{A}_\mathrm{    ref}$, we need
    to show that every model of $\Sigma_\mathrm{ref}$ lies in
    $\mathscr{A}_\mathrm{ref}$.  Before we can do this, we need some
    further consequences of $\Sigma_\mathrm{ref}$.

In the identities that occur in the next lemma we use $\overline{u}$,
where $u$ is either $x$ or $xyx$, to denote either $u$ or $u'$.
We assume that the meaning of the operation $\overline{\phantom{u}}$
is fixed within each identity: either it changes nothing or it adds $'$
to all its arguments.
    \begin{lem}\label{lem:Sigmacons}
    The following identities all follow from
    $\Sigma_\mathrm{ref}$\up:
    \begin{itemize}
        \item $(\overline{x}u_{1})'u_{2}xyx\approx
        (\overline{xyx}u_{1})'u_{2}xyx$\up;
        \item $xyxu_{2}(u_{1}\overline{x})'\approx
        xyxu_{2}(u_{1}\overline{xyx})'$\up;
        \item $(u_{1}\overline{x})'u_{2}xyx\approx
        (u_{1}\overline{xyx})'u_{2}xyx$\up;
        \item $xyxu_{2}(\overline{x}u_{1})'\approx
        xyxu_{2}(\overline{xyx}u_{1})'$\up,
    \end{itemize}
        where $u_{1}$ and $u_{2}$ are possibly empty unary
        semigroup words.
    \end{lem}
    \begin{proof}
    In each of the eight cases, if $u_{1}$ is empty, then the
    identity is equivalent modulo $x''\approx x$ to one in $\Sigma_\mathrm{ref}$ up to a
    change of letter names.  So we assume that $u_{1}$ is
    non-empty.  We can ensure that $u_{2}$ is
    non-empty by rewriting $u_{2}xyx$ and $xyxu_{2}$ as
    $(u_{2}xx')xyx$ and $xyx(x'xu_{2})$ respectively (a
    process we reverse at the end of each deduction).
    For the first identity
    we have $\Psi$ implies $(\overline{x}u_{1})'u_{2}xyx\sapprox
    ((\overline{x}'u_{2}xyx)'u_{1})'$, and then we use \eqref{eq:m} or
    \eqref{eq:r} to replace
    $\overline{x}$ by $\overline{xyx}$.  Reversing the application of $\SAAR$, we obtain the corresponding right
    hand side.

    The second identity is just a dual to the first so follows by
    symmetry.  Similarly, the fourth will follow from the third
    by symmetry.

    For the third identity,  Lemma \ref{lem:easycons} can be
    applied to the left hand side to get
    $\overline{x}'\overline{x}(u_{1}\overline{x})'u_{2}xyx$.  Now, the subword
    $\overline{x}'\overline{x}$ is either $x'x$ or $xx'$.  We
    will write it as
    $t(x,x')$ (where $t(x,y)$ is one of the words $xy$ or $yx$).  Using \eqref{eq:m},
    we have $t(x,x')(u_{1}\bar{x})'u_{2}xyx\approx t(xyx,x')(u_{1}\bar{x})'u_{2}xyx$.
    But the subword $t(xyx,x')(u_{1}\bar{x})'$ is of the form
    required to apply the
    second identity in the lemma we are proving.  Since this
    second
    identity has been established, we can use it to deduce
    $t(xyx,x')(u_{1}\overline{xyx})'u_{2}xyx$ and then reverse
    the procedure to get
    \[
    t(xyx,x')(u_{1}\overline{xyx})'u_{2}xyx\approx
    \overline{x}'\overline{x}(u_{1}\overline{xyx})'u_{2}xyx\approx
    (u_{1}\overline{xyx})'u_{2}xyx
    \]
    (the last equality requires a
    few extra easy steps in the $\overline{u}=u'$ case).
    \end{proof}

Recall that a unary polynomial $p(x)$ on an algebra $S$ is a
function $S\to S$ defined for each $a\in S$ by $p(a)=t(a,a_{1},\ldots,a_{n})$ where
$t(x,x_{1},\ldots,x_{n})$ is a term, and $a_{1},\ldots,a_{n}$ are
elements of $S$.  We let $P_{x}$ denote the set of all unary
polynomials on $S$.
The \emph{syntactic congruence} $\Syn(\theta)$ of an equivalence $\theta$ on
$S$ is defined to be
\[
\Syn(\theta):=\{(a,b)\mid p(a)\mathbin{\theta}p(b)\mbox{ for all }p(x)\in P_{x}\}.
\]
$\Syn(\theta)$ is known to be the largest congruence of $S$
contained in $\theta$ (see \cite{alm} or \cite{CDFJ}). It is very
well known that for standard semigroups, one only needs to
consider polynomials $p(x)$ built from the semigroup words
$x,x_{1}x,xx_{1},x_{1}xx_{2}$ (see \cite{how} for example).  In
fact there is a similar---though more complicated---reduction for
the variety defined by $\Sigma_\mathrm{ref}$ (and more generally
still $\Psi$). This can be gleaned fairly easily from Proposition
\ref{pro:nf} (see \cite{CDFJ} for a general approach for
establishing this), however we do not need an explicit formulation
of it here, and so omit any proof.

    We may now prove the key lemma, a variation of  \cite[Lemma 3.2]{HKMST}.
    \begin{lem}\label{lem:kub}
    Every model of $\Sigma_\mathrm{ref}$ \up(within the  variety of unary semigroups\up)
    is a subdirect product of members of $\A(\mathsf{G}_\mathrm{ref})\cup\mathscr{SB}$.
    \end{lem}
    \begin{proof}
    Let ${S}\models \Sigma_\mathrm{ref}$.  If $S$
    is the one element semigroup we are done.  Now assume that
    $|S|>1$.
    We need to show that for every pair of distinct elements
    $a,b\in S$ there is a homomorphism from $S$ onto a
    square band or an adjacency semigroup $\A(\Gg)$ for some
    $\Gg\in\underline{\mathsf{G}}_\mathrm{ref}$.

    For each element $z\in S$, we let $I_{z}:=\{u\in S\mid
    z\notin S^{1}uS^{1}\}$, in other words, $I_{z}$ is the ideal
    consisting of all elements that do not divide $z$.  Note that
    $I_{z}$ is closed under the reversion operation (since $u'$
    divides $u$).  Define equivalence relations $\rho_{z}$ and
    $\lambda_{z}$ on ${S}$:
    \[
    \rho_{z}:=\{(x,y)\in S\times S\mid (\forall t\in SzS)\quad
    xt\equiv yt \mod I_{z}\};
    \]
    \[
    \lambda_{z}:=\{(x,y)\in S\times S\mid (\forall t\in SzS)\quad
    tx\equiv ty \mod I_{z}\}.
    \]
    So far the proof is identical to that of \cite[Lemma
    3.2]{HKMST}.  In the semigroup setting, both $\rho_{z}$ and
    $\lambda_{z}$ are congruences, however this is no longer true
    in the unary semigroup setting.  Instead, we replace $\rho_{z}$ and
    $\lambda_{z}$ by their syntactic congruences $\Syn(\rho_{z})$
    and $\Syn(\lambda_{z})$.

    Let $a$ and $b$ be distinct elements of $S$.  Our goal is to
    show that one of the congruences $\Syn(\rho_{a})$,
    $\Syn(\rho_{b})$,  $\Syn(\lambda_{a})$ and $\Syn(\lambda_{b})$ separate $a$ and $b$, and
    that ${S}/\Syn(\rho_{z})$ and ${S}/\Syn(\lambda_{z})$ are
    isomorphic to a square band or an adjacency semigroup of a
    reflexive graph.  The first part is essentially identical to
    a corresponding part of the proof of
    \cite[Lemma 3.2]{HKMST}.  We include it for completeness only.

    First suppose that $a\notin SbS$.  So $b\in I_{a}$.  Choose
    $t=a'a\in SaS$ so that $a=at\not\equiv bt\mod I_{a}$.  Hence
    $(a,b)\notin \hat\rho_{a}$.  Now suppose that $SaS=SbS$, so that
    $a$ and $b$ lie in the same $\mathrsfs{J}$-class
    $SaS\backslash I_{a}$ of ${\bf S}$.  One of the following
    two equalities must fail: $ab'b=b$ or $aa'b=a$ for otherwise
    $a=aa'b=aa'ab'b=ab'b=b$.  Hence as neither $a$ nor $b$ is in
    $I_{a}=I_{b}$, we have either
    $(a,b)\notin\rho_{a}\supseteq \hat\rho_{a}$ or
    $(a,b)\notin\lambda_{a}\supseteq \hat\lambda_{a}$.

    Now it remains to prove that ${S}/\Syn(\rho_{z})$ and
    ${S}/\Syn(\lambda_{z})$ are adjacency semigroups or square bands.
    Lemmas~\ref{lem:reflexive}, \ref{lem:reflexive2} and
    \ref{lem:easycons} show that it
    suffices to prove that the underlying semigroup of ${S}/\Syn(\rho_{z})$ is completely
    0-simple or completely simple.  We look at the $\Syn(\rho_{z})$ case only (the
    $\Syn(\lambda_{z})$ case follows by symmetry).  Now it does no harm to assume
    that $I_{z}$ is empty or $\{0\}$, since $v,w\in I_{z}$
    obviously implies that $(v,w)\in \Syn(\rho_{z})$.  Hence $K_{z}:=SzS/(I_{z}\cap
    SzS)$ is a 0-simple semigroup or a simple semigroup.  Since $S$ is periodic
    (by identity \eqref{eq:aperiodic} of $\Sigma_\mathrm{ref}$), we have that $K_{z}$ is completely
    0-simple or completely simple.  We need to prove that every element of
    $S\backslash I_{z}$ is  $\Syn(\rho_{z})$-related to a member of $SzS\backslash I_{z}$.

    Let $c\in S$.  If $c\in SzS$ or $c\in I_{z}$ we are done, so let us assume
    that $c\notin SzS\cup I_{z}$.  So $z=pcq$ for some $p,q\in
    S^{1}$.  So $z=pcqz'pcq$.  Put $w=qz'p$.  Note that $w\in
    SzS$ and $cwc\neq 0$.  0ur goal is to show that
    $c\mathbin{\Syn(\rho_{z})}cwc$.  Let $s(x,\vec{y})$ be
    any unary semigroup word in some variables $x,y_{1},\ldots$ and
    let $t\in SaS$.  We need to prove that for any
    $\vec{d}$ in $S^{1}$ we have
    $s(c,\vec{d}\,)t\equiv s(cwc,\vec{d}\,)t$
    modulo $I_{z}$.  Write $t$ as $ucwcv$, which is possible
    since both $t$ and $cwc$ are $\mathrsfs{J}$-related.  (Note
    that modulo the identity $xx'x\approx x$ we may assume both
    $u$ and $v$ are nonempty.)  We want
    to obtain
    \begin{equation}
        s(c,\vec{d}\,)ucwcv=s(cwc,\vec{d}\,)ucwcv.\label{eq:kub}
    \end{equation}

    Now using Corollary \ref{cor:refnf}, we may rewrite
    $s(c,\vec{d}\,)$ as a word in which each application of $'$
    covers either a single variable or a word of the form $gh$
    where $g,h$ are either letters or $'$ applied to a letter.
    There may be many occurrences of $c$ in this word.  We show
    how to replace an arbitrary one of these by $cwc$ and by
    repeating this to each of these occurrences we will achieve the desired
    equality \eqref{eq:kub}.  Let us fix some occurrence of $c$.
    So we may consider the expression $s(c,\vec{d}\,)ucwc$ as
    being of one of the following forms:
    $w_{1}cw_{2}cwc$; $w_{1}c'w_{2}cwc$; $w_{1}(cz)'w_{2}cwc$;
    $w_{1}(c'z)'w_{2}cwc$; $w_{1}(zc')'w_{2}cwc$;
    $w_{1}(zc)'w_{2}cwc$.  In each case, we can make the required
    replacement using a single  application of Lemma
    \ref{lem:Sigmacons}.  This gives equality \eqref{eq:kub}, which
    completes the proof.
    \end{proof}

    As an immediate corollary we obtain the following result.
    \begin{cor}\label{cor:basis for ref}
    The identities $\Sigma_\mathrm{ref}$ are an identity basis for
    $\mathscr{A}_\mathrm{ref}$.
    \end{cor}

    Let $\mathscr{SL}$ denote is the variety  generated by the adjacency semigroup
    over the 1-vertex looped graph and let $\mathscr{U}$ denote the variety generated
    by adjacency semigroups over single block equivalence relations (equivalently,
    $\mathscr{U}$ is the variety  generated by the adjacency semigroup over the universal
    relation on a 2-element set).  Recall that $\mathscr{SB}$ denotes the
    variety of square bands.

    \begin{lem}\label{lem:Sjoin}
    $\mathscr{SL}\vee\mathscr{SB}= \mathscr{U}$.
    \end{lem}
    \begin{proof}
        The direct product of the semigroup $\A(\underline{\mathbf{1}})$ with
        an $I\times I$ square band has a unique maximal ideal
        and the corresponding Rees quotient is (isomorphic to) the adjacency
        semigroup over the universal relation on $I$.  So
        $\mathscr{SL}\vee\mathscr{SB}\supseteq \mathscr{U}$.
        However
        if $|I|\geq 2$, and $\mathsf{U}_{I}$ denotes the
        universal relation on $I$, then the adjacency semigroup $\A(\mathsf{U}_{I})\in
        \mathscr{U}$
        contains as subalgebras both $\A(\underline{\mathbf{1}})$ (a generator
        for $\mathscr{SL}$) and the $I\times I$ square band (a generator for
        $\mathscr{SB}$).  So $\mathscr{SL}\vee\mathscr{SB}\subseteq \mathscr{U}$.
    \end{proof}

    \begin{lem}\label{lem:squareband}
    Let $\mathscr{V}$ be a subvariety of $\mathscr{A}_\mathrm{ref}$
    containing the variety $\mathscr{SB}$.  Either $\mathscr{V}=\mathscr{SB}$ or
    $\mathscr{V}\supseteq \mathscr{U}$ and
    $\mathscr{V}=\mathbb{HSP}(\A(K))$ for some class of
    \up(necessarily reflexive\up) graphs $K$.
    \end{lem}
    \begin{proof}
        Let ${\bf A}$ be a nonfinitely generated
        $\mathscr{V}$-free algebra.  If ${\bf A}\models xyx\approx
        x$ then $\mathscr{V}$ is equal to $\mathscr{SB}$.  Now say
        that $xyx\approx x$ fails on ${\bf A}$.
        Lemma \ref{lem:kub}, shows that ${\bf A}$ is a subdirect
        product of some family $J$ of adjacency semigroups
        and square bands.  Note that we have
        $\mathscr{V}=\mathbb{HSP}(J)$.  Our
        goal is to replace all square bands in $J$ by
        adjacency semigroups over universal relations.

        Since $xyx\approx
        x$ fails on ${\bf A}$, at least one of the subdirect factors of
        ${\bf A}$ is an adjacency semigroup that is not the one
        element algebra.  Hence $\mathscr{V}$ contains the semigroup
        $\A(\underline{\bf 1})$.  By Lemma~\ref{lem:Sjoin},
        $\mathscr{V}$ contains $\mathscr{U}$.  Now replace
        all square bands in $J$ by the adjacency semigroup of
        a universal relation of some set of size at least 2, and denote
        the corresponding class by $\bar{J}$; let
        $\underline{\mathsf{G}}_{\bar{J}}$ denote the
        corresponding class of graphs.  Then
        $\mathscr{V}=\mathbb{HSP}(J)
        =\mathbb{HSP}(\bar{J})=\mathbb{HSP}(\A(\underline{\mathsf{G}}_{\bar{J}}))$.
    \end{proof}

    Now we may complete the proof of Theorem \ref{thm:fundref}.
    \begin{proof}
    Theorem \ref{thm:fund} shows the map
    $\iota$ described in Theorem \ref{thm:fund} is an order
    preserving injection from $L(\underline{\mathsf{G}}_\mathrm{ref})$
    to $L(\mathscr{A}_\mathrm{ref})$.  Now we show that it is a
    surjection.  That is, every subvariety of $\mathscr{A}_\mathrm{ref}$
    other than $\mathscr{SB}$ is the image  under $\iota$ of some uH class
    of reflexive graphs.  Lemma~\ref{lem:squareband} shows this is true if
    $\mathscr{SB}\subseteq \mathscr{V}$.  However, if the square
    bands in $\mathscr{V}$ are all trivial, then Lemma
    \ref{lem:kub} shows that either $\mathscr{V}$ is the
    trivial variety (and equal to $\iota(\{{\bf 0}\})$) or
    $\omega$-generated $\mathscr{V}$-free algebra is a
    subdirect product of members of
    $\A(\underline{\mathsf{G}}_\mathrm{ref})$.  Let $F$
    be a set consisting of the subdirect factors and
    $\underline{\mathsf{G}}_{F}$ the
    corresponding graphs.  Then
    $\mathscr{V}=\mathbb{HSP}(F)=\iota(\mathbb{ISP^{+}P}_\mathrm{    u}(\underline{\mathsf{G}}_{F}))$.  To show that $\iota$ is a
    lattice isomorphism, it will suffice to show that $\iota$
    preserves joins, since meets follow from the fact that
    $\iota$ is an order preserving bijection.

    Let $\bigvee_{i\in I}\mathsf{\underline{R}}_{i}$ be some join in
    $L^{+}$.  First assume that $\mathsf{\underline{S}}$ is not amongst the
    $\mathsf{\underline{R}}_{i}$.  Then
    \begin{multline*}
    \mathbb{HSP}(\A(\bigvee_{i\in
    I}\mathsf{R}_{i}))=\mathbb{HSP}(\A(\mathbb{ISP^{+}P}_\mathrm{u}(\bigcup_{i\in
    I}\mathsf{R}_{i})))\\
    =\mathbb{HSP}(\mathbb{HSP}(\bigcup_{i\in
    I}\A(\mathsf{R}_{i})))=
    \bigvee_{i\in I}\mathbb{HSP}(\A(\mathsf{R}_{i})).
    \end{multline*}
    If $\mathsf{\underline{S}}$ is amongst the
    $\mathsf{\underline{R}}_{i}$ then either the join is a join of
    $\mathsf{\underline{S}}$ with the trivial uH class $\{{\bf
    0}\}$ (and the join is obviously preserved by $\iota$), or using Lemma
    \ref{lem:Sjoin}, we can replace $\mathsf{\underline{S}}$ by
    the uH class of universal relations, and proceed as above.
    This completes the characterization of
    $\mathscr{L}(\mathscr{A}_\mathrm{ref})$.

    Next we must show that a
    class $K$ of graphs generates a finitely
    axiomatizable uH class if and only if
    $\mathbb{HSP}(\A(K))$ is a finitely axiomatizable.
    The ``only if'' case is Corollary \ref{cor:NFB}.  Now say that
    $K$ has a finite basis for its uH sentences.
    Following the methods of Subsection \ref{subsec:fund}, we may
    construct a finite set $\Xi$ of identities such that an adjacency
    semigroup $\A$ lies in $\mathbb{HSP}(\A(K))$ if
    and only if $\A\models \Xi$.  We claim that $\Sigma_\mathrm{ref}\cup\Xi$
    is an identity basis for
    $\mathbb{HSP}(\A(K))$.  Indeed, if $S$ is a unary semigroup
    satisfying $\Sigma_\mathrm{ref}\cup\Xi$, then by Lemma~\ref{lem:kub},
    $S$ is a subdirect product of adjacency  semigroups (or possibly square
    bands) satisfying $\Xi$.  So these adjacency semigroups lie in $\mathbb{HSP}(\A(K))$,
    whence so does ${S}$.

    The proof that $\iota$ preserves the property of being
    finitely generated (and being nonfinitely generated) is very
    similar and left to the reader.
    \end{proof}

\section{Applications}
\label{sec:applications}
The universal Horn theory of graphs is reasonably well developed, and the link
to unary Rees matrix semigroups that we have just established provides numerous
corollaries. We restrict ourselves to just a few ones which all are based on
the examples of uH classes presented in Section~\ref{sec:graphs}.

We start with presenting finite generators for unary semigroup varieties
that we have considered.
    \begin{pro}
        The varieties $\mathscr{A}$, $\mathscr{A}_\mathrm{symm}$\up, and
        $\mathscr{A}_\mathrm{ref}$ are generated by
        $\A(\mathsf{G}_{1})$\up,  $\A(\mathsf{S}_{1})$ and
        $\A(\mathsf{R}_{1})$ respectively.
    \end{pro}
    \begin{proof}
        This follows from Theorem~\ref{thm:fund} and Examples~\ref{eg:graphgen}, \ref{eg:symmetric}, and~\ref{eg:refgen}.
    \end{proof}

Observe that the generators are of fairly modest size, with 17, 17 and
10~elements respectively.

Recall that $\Cg_{3}$ is a 5-vertex graph generating the uH class of
all $3$-colorable graphs (Example \ref{eg:k-color}, see also Fig.\,\ref{pic:c2 and c3}).
    \begin{pro}\label{eg:C3}
        The finite membership problem for the variety generated by the
        $26$-element unary semigroup $\A(\Cg_{3})$ is \textsf{NP}-hard.
    \end{pro}
    \begin{proof}
        By Theorem \ref{thm:fund} $\A(\Gg)\in\mathbb{HSP}_\mathrm{fin}(\A(\Cg_{3}))$
        if and only if $\Gg$ is $3$-colorable, a known \textsf{NP}-complete problem, see
        \cite{garey}. Of course, the construction of $\A(\Gg)$ can be made in
        polynomial time, so this is a polynomial reduction.
    \end{proof}
A similar (but more complicated) example in the plain semigroup setting has been
found in~\cite{jacmck}. Observe that we do not claim that the finite membership problem
for $\mathbb{HSP}(\A(\Cg_{3}))$ is \textsf{NP}-complete since it is not clear whether
or not the problem is in \textsf{NP}.

One can also show that the equational theory of $\A(\Cg_{3})$ is co-\textsf{NP}-complete.
(It means that the problem whose instance is a unary semigroup identity $u\approx v$
and whose question is whether or not $u\approx v$ holds in $\A(\Cg_{3})$ is
co-\textsf{NP}-complete.) This follows from the construction of identities modelling
uH sentences in Subsection~\ref{subsec:fund}. The argument is an exact parallel to that
associated with \cite[Corollary 3.8]{jacmck} and we omit the details.

    \begin{pro}\label{cor:NFB}
        If $K$ is a class of graphs without a finite basis of
        uH sentences, then $\A(K)$ is without a finite basis of identities.
        If $K$ is a class of graphs whose uH class has \up(infinitely many\up)
        uncountably many sub-uH classes, then the variety generated by
        $\A(K)$ has \up(infinitely many\up) uncountably many
        subvarieties.
    \end{pro}

    \begin{proof}
        This is an immediate consequence of Theorem \ref{thm:fund}.
    \end{proof}

In particular,  recall the 2-vertex graph $\Sg_{2}$ of Example \ref{eg:simplegraphgen},
and let $\Kg_{2}$ denote the 2-vertex complete simple graph.
    \begin{cor}\label{eg:uncountable}
        There are uncountably many varieties between the variety
        generated by $\A(\Sg_{2})$ and that generated by $\A(\Kg_{2})$.
        The statement is also true if $\Sg_{2}$ is replaced by any simple graph
        that is not $2$-colorable.
     \end{cor}

    \begin{proof}
        The first statement follows from Theorem \ref{thm:fund}
        and statements in Example~\ref{eg:simplegraphgen}.  The
        second statement follows similarly from statements in
        Example~\ref{eg:simple}.
    \end{proof}

    Note that the underlying semigroup of $\A(\Sg_{2})$ is simply the familiar
    semigroup ${A}_{2}$, see Subsection~\ref{subsec:equational}. The subvariety
    lattice of the semigroup variety generated by $A_{2}$ is reasonably well
    understood (see Lee and  Volkov~\cite{leevol}).  This variety contains all
    semigroup varieties generated by completely 0-simple semigroups with
    trivial subgroups but has only countably many subvarieties, all of which
    are finitely axiomatized (see Lee~\cite{lee}).

    Theorem \ref{thm:fundref} reduces the study of the subvarieties
    of $\mathscr{A}_\mathrm{ref}$ to the study of uH classes of reflexive
    graphs.  This class of graphs does not seem to have been as
    heavily investigated as the antireflexive graphs, but contains
    some interesting examples.

Recently Trotta~\cite{trotta} has disproved a claim made in~\cite{sizyi}
by showing that there are uncountably many uH classes of reflexive
antisymmetric graphs. From this and Theorem \ref{thm:fundref} we
immediately deduce:

    \begin{pro}\label{eg:infinite}
    The unary semigroup variety $\mathscr{A}_\mathrm{ref}$ has uncountably many
    subvarieties.
    \end{pro}

\begin{figure}[th]
\begin{center}
    \begin{pspicture}(8,4)(-2,0)
      \cnodeput(0,0){a}{}
      \cnodeput(1,1){b}{}
      \cnodeput(0,2){c}{}
      \cnodeput(2,2){d}{}
      \cnodeput(1,3){e}{}
      \cnodeput(0,4){h}{}
      \psset{arrows=-,nodesep=0}
      \ncline{a}{b}
      \ncline{b}{c}
      \ncline{b}{d}
      \ncline{c}{e}
      \ncline{e}{h}
      \ncline{d}{e}
      \put(0.3,-0.1){$\{\underline{\bf 0}\}$}
      \put(1.3,0.9){$\mathbb{I}(\{\underline{\bf 1}\})$}
      \put(-3.2,1.9){Single block}
      \put(-3.2,1.5){equivalence relations}
      \put(2.3,1.9){Antichains}
      \put(-3,2.9){Equivalence relations}
      \put(-3.2,3.9){All reflexive and}
      \put(-3.2,3.5){symmetric graphs}
      \cnodeput(6,0){a1}{}
      \cnodeput(7,1){b1}{}
      \cnodeput(6,2){c1}{}
      \cnodeput(8,2){d1}{}
      \cnodeput(7,3){e1}{}
      \cnodeput(6,4){h1}{}
      \cnodeput(5,1){s}{}
      \psset{arrows=-,nodesep=0}
      \ncline{a1}{s}
      \ncline{s}{c1}
      \ncline{a1}{b1}
      \ncline{b1}{c1}
      \ncline{b1}{d1}
      \ncline{c1}{e1}
      \ncline{e1}{h1}
      \ncline{d1}{e1}
      \put(6.3,-0.1){$\mathscr{T}$}
      \put(7.3,0.9){$\mathscr{SL}$}
      \put(8.3,1.9){$\mathscr{BR}$}
      \put(4.5,2.2){$\mathscr{SL}\vee\mathscr{SB}$}
      \put(7.3,2.9){$\mathscr{BR}\vee\mathscr{SB}$}
      \put(5.3,0.9){$\mathscr{SB}$}
      \put(6.3,3.9){$\mathscr{CSR}$}
    \end{pspicture}
    \caption{The lattice of uH classes of reflexive symmetric graphs vs
    the lattice of varieties of strict regular semigroups}\label{pic:refsym}
\end{center}
\end{figure}
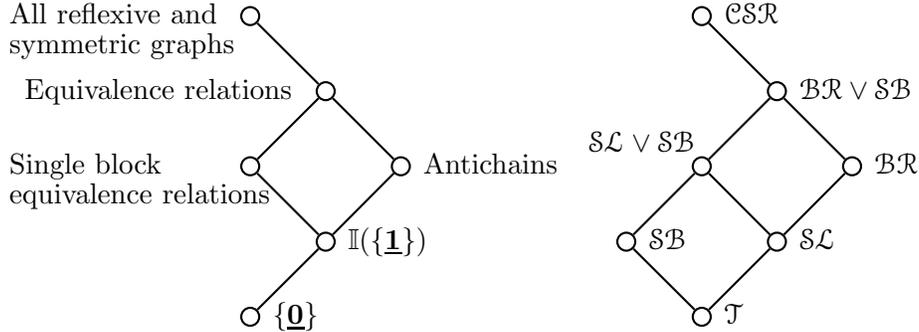

In contrast, it is easy to check that there are only 6 uH classes of
reflexive symmetric graphs, see~\cite{ben} for example. The lattice they
form is shown in Fig.\,\ref{pic:refsym} on the left.
Theorem~\ref{thm:fundref} then implies that the subvariety lattice
of the corresponding variety of unary semigroups contains 7 elements
(it is one of the cases when the ''extra'' variety $\mathscr{SB}$ of
square bands comes into the play); the lattice is shown in Fig.\,\ref{pic:refsym}
on the right. The variety is generated by the adjacency semigroup of the graph
$\mathsf{RS}_{1}$ of Example \ref{eg:refgen} and is nothing but the variety
$\mathscr{CSR}$ of so-called combinatorial strict regular *-semigroups which
have been one of the central objects of study in \cite{aui}. The other
join-indecomposable varieties in Fig.\,\ref{pic:refsym} are the trivial variety
$\mathscr{T}$, the variety $\mathscr{SL}$ of semilattices with identity map
as the unary operation, and the variety $\mathscr{BR}$ of combinatorial strict
inverse semigroups.

The main results of \cite{aui} consisted in providing a finite identity basis
for $\mathscr{CSR}$ and determining its subvariety lattice. We see that the
latter result is an immediate consequence of Theorem~\ref{thm:fundref}.
A finite identity basis for $\mathscr{CSR}$ can be obtained by adding the
involution identity \eqref{eq:involution} to the identity basis $\Sigma_\mathrm{ref}$
of the variety $\mathscr{A}_\mathrm{ref}$, see Corollary~\ref{cor:basis for ref}.
(The basis constructed this way is different from that given in \cite{aui}.)

    \begin{eg}\label{eg:po}
    The adjacency semigroup $\A(\underline{\bf 2})$ of the two element
    chain $\underline{\bf 2}$ \up(as a partial order\up) generates
    a variety with a lattice of  subvarieties isomorphic to the four
    element chain.  The  variety is a cover of the variety $\mathscr{BR}$
    of combinatorial strict inverse semigroups.
    \end{eg}

    \begin{proof}
    This follows from  Example \ref{eg:preorders}, Theorem~\ref{thm:fundref}
    and the fact that the uH class of universal relations is not a sub-uH class
    of the partial orders (so $\mathscr{SB}$ is not a subvariety of
    $\mathbb{HSP}(\A(\underline{\bf 2}))$.
    \end{proof}

    The underlying semigroup of $\A(\underline{\bf 2})$ is again the
    semigroup ${A}_{2}$. Thus, Example \ref{eg:po} makes an
    interesting contrast to Corollary~\ref{eg:uncountable}.

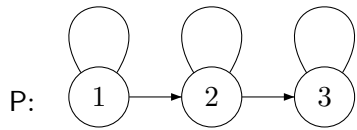
\begin{figure}[th]
\begin{center}
\begin{picture}(60,20)
\unitlength=.6mm
\node(B)(0,0){1}
\node(C)(25,0){2}
\node(D)(50,0){3}
\drawedge(B,C){}
\drawedge(C,D){}
\put(-20,-3){$\Pg$:}
\gasset{AHnb=0}
\drawloop(B){}
\drawloop(C){}
\drawloop(D){}
\end{picture}
    \caption{Generator for a limit uH class}\label{pic:limit}
\end{center}
\end{figure}

 For our final application, consider the $3$-vertex graph $\Pg$ shown in Fig.\,\ref{pic:limit}.

It is known (see \cite{ben}) and easy to verify that the uH-class $\mathbb{ISP^{+}P}_\mathrm{u}(\Pg)$
is not finitely axiomatizable and the class of partial orders is the unique maximal sub-uH class of
$\mathbb{ISP^{+}P}_\mathrm{u}(\Pg)$. Recall that a variety $\mathscr{V}$ is said to be a \emph{limit}
variety if $\mathscr{V}$ has no finite identity basis while each of its proper subvarieties is finitely
based. The existence of limit varieties is an easy consequence of Zorn's lemma but concrete examples
of such varieties are quite rare. We can use the just registered properties of the graph $\Pg$ in order
to produce a new example of a finitely generated limit variety of I-semigroups.

    \begin{pro}\label{pro:Plimit}
    The variety $\mathbb{HSP}(\A(\Pg))$ is a limit variety whose subvariety lattice is
    a $5$-element chain.
    \end{pro}

    \begin{proof}
    This follows from Theorem \ref{thm:fundref} and Example \ref{eg:preorders}.
    \end{proof}

\section*{Conclusion}

We have found a transparent translation of facets of universal Horn logic
into the apparently much more restricted world of equational logic.
A general translation of this sort has been  established for uH classes
of arbitrary structures (even partial structures) by the first author~\cite{jacFlat}.
We think however that the special case considered in this paper is of interest
because it deals with very natural objects on both universal Horn logic and
equational logic sides.

We have shown that the unary semigroup variety $\mathscr{A}_\mathrm{ref}$ whose
equational logic captures the universal Horn logic of the reflexive graphs is
finitely axiomatizable. The question of whether or not the same is true for
the variety $\mathscr{A}$ corresponding to all graphs still remains open.
A natural candidate for a finite identity basis of $\mathscr{A}$ is the
system consisting of the identities $(\Psi)$ and \eqref{eq:fix}--\eqref{eq:l},
see Section~\ref{sec:regular}.

\bibliographystyle{amsplain}

\end{document}